\documentclass[11pt,final]{amsart}
\usepackage[text={5in,8in},centering]{geometry}
\usepackage{graphicx}
\usepackage{amssymb, amsmath, amsthm}
\usepackage[colorlinks=true, citecolor=black, linkcolor=black, urlcolor=black]{hyperref}
\usepackage{mathrsfs}  

\usepackage{bm}

\usepackage{ifdraft}
\ifoptionfinal{
\usepackage[disable]{todonotes}
}{
\usepackage[norefs, nocites]{refcheck}

\usepackage[notref, notcite]{showkeys}
\usepackage[bordercolor=white, color=white]{todonotes}
}
\newcommand{\HOX}[1]{\todo[noline, size=\footnotesize]{#1}}

\makeatletter\providecommand\@dotsep{5}\def\listtodoname{List of Todos}\def\listoftodos{\hypersetup{linkcolor=black}\@starttoc{tdo}\listtodoname\hypersetup{linkcolor=blue}}\makeatother

\usepackage{hyperref}
\usepackage{etoolbox}

\newtheorem{theorem}{Theorem}[section]
\newtheorem{lemma}[theorem]{Lemma}
\newtheorem{corollary}[theorem]{Corollary}
\newtheorem{proposition}[theorem]{Proposition}
\newtheorem{definition}[theorem]{Definition}

\theoremstyle{remark}
\newtheorem{remark}{Remark}

\numberwithin{equation}{section}

\def\R{\mathbb R}

\def\N{\mathbb N}

\renewcommand{\leq}{\leqslant}
\renewcommand{\geq}{\geqslant}

\def\linspan{\mathrm{span}}

\def\p{\partial}
\DeclareMathOperator{\supp}{supp}
\newcommand{\pair}[1]{\left\langle #1 \right\rangle}

\DeclareMathOperator{\spn}{span}
\newcommand\KN{~\wedge\!\!\!\!\!\!\!\!\;\bigcirc~}
\let\div\relax
\DeclareMathOperator{\div}{div}

\DeclareMathOperator{\Hess}{Hess}

\newcommand*\xbar[1]{%
   \hbox{%
     \vbox{%
       \hrule height 0.5pt 
       \kern0.5ex
       \hbox{%
         \ensuremath{#1}%
       }%
     }%
   }%
} 

\title[Lorentzian Calder\'{o}n problem under curvature bounds]{Lorentzian Calder\'{o}n problem under curvature bounds}
\author[Alexakis]{Spyros Alexakis}
\address{Department of Mathematics, University of Toronto, 40 St George St, Toronto, ON, Canada M5S 2E4.}
\email{alexakis@math.toronto.edu}
\author[Feizmohammadi]{Ali Feizmohammadi}
\address{Department of Mathematics, University College London, 
Gower Street, London UK, WC1E 6BT.}
\email{a.feizmohammadi@ucl.ac.uk}
\author[Oksanen]{Lauri Oksanen}
\address{Department of Mathematics, University College London, 
Gower Street, London UK, WC1E 6BT.}
\email{l.oksanen@ucl.ac.uk}


\begin{document}

\maketitle
\begin{abstract}
We introduce a method of solving inverse boundary value problems for wave equations on Lorentzian manifolds,
and show that zeroth order coefficients can be recovered 
under certain curvature bounds. 
The set of Lorentzian metrics satisfying the curvature bounds has a non-empty interior in the sense of arbitrary, smooth perturbations of the metric, whereas all previous results on this problem impose conditions on the metric that force it to be real analytic with respect to a suitably defined time variable. The analogous problem on Riemannian manifolds is called the Calder\'on problem, and in this case the known results require the metric to be independent of one of the variables. 
Our approach is based on a new unique continuation result in the exterior of the double null cone emanating from a point.
The approach shares features with the classical Boundary Control method, and can be viewed as a generalization of this method to cases where no real analyticity is assumed. 
\end{abstract}

\section{Introduction}

Let $(\mathcal M,\textsl{g})$ be a connected, smooth Lorentzian manifold with timelike boundary. We write $1+n$ for the dimension and $\mathcal M^{\textrm{int}}$, $\p \mathcal M$ for the interior and boundary of $\mathcal M=\mathcal M^{\textrm{int}}\cup \p\mathcal M$. Let $V \in C^\infty(\mathcal M)$ and consider the Cauchy data set 
    \begin{align}
    \label{C_V}
\mathscr C(V) = \{(u, \p_\nu u)|_{\p \mathcal M} : 
\text{$u \in C^\infty(\mathcal M)$ and $\Box u + V u = 0$ on $\mathcal M$}\},
    \end{align}
where $\Box$ is the canonical wave operator on $(\mathcal M,\textsl{g})$ and $\nu$ is the exterior unit normal vector field on $\p \mathcal M$.
We call the problem to find $V$ given $\mathscr C(V)$
the Lorentzian Calder\'{o}n problem. The classical Calder\'{o}n problem has the same formulation except that $(\mathcal M,\textsl{g})$ is a smooth Riemannian manifold with boundary and $\Box$ is replaced by the Laplacian on $(\mathcal M,\textsl{g})$.

Another version of the Calder\'{o}n problem is to find
$\textsl{g}$ up to an isometry given the Cauchy data set 
    \begin{align*}
\{(u, \p_\nu u)|_{\p \mathcal M} : 
\text{$u \in C^\infty(\mathcal M)$ and $\Box u = 0$ on $\mathcal M$}\}.
    \end{align*}
If $n>1$, then in a fixed conformal class, the latter problem reduces to the former one by using a gauge transformation.
Writing $\Box_c$ for the wave operator with respect to the conformally scaled metric $c \,\textsl{g}$,
the function
    \begin{align}\label{gauge}
w = c^{(n-1)/4} u
    \end{align}
satisfies the equation $\Box w + V w = 0$, with $V = -c^{-(n-1)/4} \Box c^{(n-1)/4}$,
if the function $u$ satisfies the equation $\Box_{c} u = 0$. 

Without any further assumptions on $(\mathcal M,\textsl{g})$ 
the Calder\'{o}n problem is wide open, 
in both the Lorentzian and Riemannian cases, and regardless whether $\textsl{g}$ or $V$ is to be determined. 
We will now formulate the geometric assumptions under which we show that $\mathscr C(V)$ determines $V$ in the Lorentzian case.

In order to be able to solve the wave equation on $\mathcal M$, and thus guarantee that $\mathscr C(V)$ has a rich structure, we assume that 
\begin{enumerate}
\item[(H1)] There is a smooth, proper, surjective temporal function $\tau : \mathcal M \to \R$.
\end{enumerate}
By $\tau$ being a temporal function we mean that its differential $d\tau$ is timelike. Proper is used in the topological sense, that is, inverse images of compact subsets are compact under $\tau$.

Hypothesis~(H1), together with $\p \mathcal M$ being timelike,
implies that the wave equation can be solved on $\mathcal M$, see \cite[Theorem 24.1.1]{Ho3}.
This theorem does not require $\tau$ to be surjective. However, surjectivity guarantees that $\mathcal M$ is diffeomorphic to a cylinder, a natural feature in view of the previous results discussed in Section \ref{previous_lit} below. This assumption avoids some technical complications.  

When wave equations are considered on a Lorentzian manifold without boundary $(\widetilde{\mathcal M},\textsl{g})$, it is typically assumed that the manifold is globally hyperbolic. 
In this case $\widetilde{\mathcal M}$ is isometric to a cylinder $\R \times \widetilde{M_0}$
with a metric of the form
    \begin{align}\label{ortho_split}
c(t,x) (-dt^2 + g_0(t,x)),
    \end{align}
where $\widetilde{M_0}$ is a smooth manifold, $c$ is a smooth positive function, and $g_0(t,\cdot)$ is a family of smooth Riemannian metrics on $\widetilde{M_0}$ that depend smoothly on the variable $t \in \R$, see \cite{BS}.
Observe that (H1) is satisfied for $\tau = t$ on $\mathcal M = \R \times M_0$ with $M_0 \subset \widetilde{M_0}$ an open, bounded set with smooth boundary.
In general, we will show that (H1), together with $\p \mathcal M$ being timelike, implies that $\mathcal M$ is diffeomorphic to a cylinder $\R \times M_0$, where $M_0$ is a compact smooth manifold with boundary, and that over compact subsets of $\mathcal M$, the metric is isometric to a metric of the form (\ref{ortho_split}).

Let us now turn to our main assumption on the curvature of $(\mathcal M,\textsl{g})$.
We fix the signature convention $(-,+,\ldots,+)$, let $R$ stand for the curvature tensor on $(\mathcal M,\textsl{g})$, and recall the following definition from \cite{AH98}.
\begin{definition}
\label{def_K}
For $K \in \R$, we write $R \le K$ if
    \begin{align*}
\textsl{g}(R(X,Y)Y,X) \leq K \left( \textsl{g}(X,X)\textsl{g}(Y,Y)-\textsl{g}(X,Y)^2 \right)
    \end{align*}
for all $X,Y \in T_p\mathcal M$ and $p \in \mathcal M$.
\end{definition}
In the Riemannian case the curvature bound $R \le K$
is equivalent with 
$\mathrm{Sec}(X,Y) \leq K$
for all linearly independent $X,Y \in T_p\mathcal M$ and $p \in \mathcal M$,
where $\mathrm{Sec}$ is the sectional curvature,
    \begin{align*}
\mathrm{Sec}(X,Y) = \frac{\textsl{g}(R(X,Y)Y,X)}{\textsl{g}(X,X)\textsl{g}(Y,Y)-\textsl{g}(X,Y)^2}.
    \end{align*}
In the Lorentzian case, these two formulations are no longer equivalent, and the latter one leads to an uninteresting theory. If $\mathrm{Sec}(X,Y) \leq K$ whenever $\mathrm{Sec}(X,Y)$ is well-defined, then the Lorentzian manifold is of constant sectional curvature \cite{K}.

A Lorentzian manifold of constant sectional curvature satisfies $R \le K$ with $K$ the constant value of the sectional curvature, but a small perturbation of a manifold with constant curvature might not satisfy $R \le K$ for any $K \in \R$. However, there are manifolds satisfying $R \le K$ such that their small perturbations satisfy the same curvature bound, see \cite{AB08}.
These examples will be discussed in more detail below (see Section~\ref{curv_cond_sec} and Section~\ref{sec_perturb}). 
We mention that there is also an equivalent description of the curvature bound $R\le K$ in terms of local triangle comparisons of the signed lengths of
geodesics \cite{AB08}. 

In order to state our main hypothesis on $(\mathcal M,\textsl{g})$, we recall that the spatial diameter of $\mathcal M$, denoted by $\textrm{Diam}(\mathcal M)$, is the supremum of the length of inextendible spacelike geodesics on $\mathcal M$. We assume that
\begin{enumerate}
\item[(H2)] 
$R\le K$ for some $K \in \R$, and if $K>0$
then $\textrm{Diam}(\mathcal M) < \frac{\pi}{2\sqrt{K}}$.
\end{enumerate}

\begin{figure}
\includegraphics[trim={5cm 2cm 3cm 0},clip,width=0.5 \textwidth]{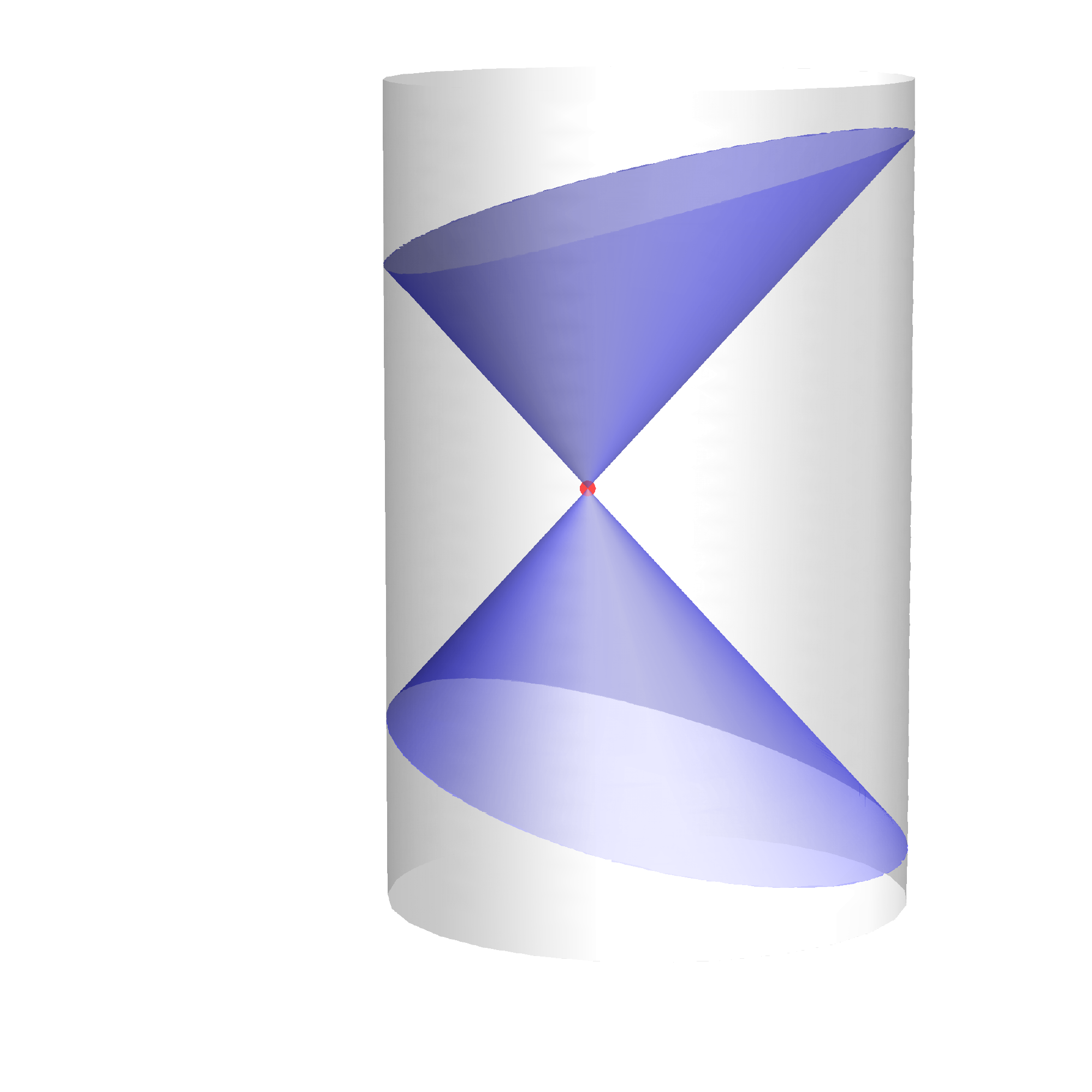}{\centering}
\caption{The schematic for the exterior of the double null cone in the setting of Minkowski geometry in $\R^{1+2}$. The point $p$ is shown in red, $\p \mathcal M$ is gray and $\p \mathscr E_p \cap \mathcal M^{\textrm{int}}$ is shown in blue.}
    \label{default}
\end{figure}

In addition to (H1)--(H2) we will make some technical assumptions. Given any $p \in \mathcal M$, let $J^+(p)$ and $J^-(p)$ denote the causal future and past of the point $p$ respectively, and define 
\begin{equation}
\label{exterior_null}
\mathscr E_p = \mathcal M \setminus (J^-(p)\cup J^+(p)).
\end{equation}
The definitions of $J^{\pm}(p)$ are recalled in Section~\ref{sec_perturb} below (see equation~\ref{causal_rel}). We call $\mathscr E_p$ the exterior of the double null cone emanating from the point $p$, see Figure~\ref{default}. We assume that
\begin{enumerate}
\item[(H3)] For any null geodesic $\gamma$ and any two points $p,q$ on $\gamma$, the only causal path between $p$ and $q$ is along $\gamma$.
For all $p \in \mathcal M$, the exponential map $\exp_p$
is a diffeomorphism from the spacelike vectors (in its maximal domain of definition) onto $\mathscr E_p$.
\end{enumerate}
We remark that (H2) has certain common features with (H3). For instance, following the proof of \cite[Proposition 11.13]{Beem}, we see that $R\le 0$ implies that null geodesics do not have any conjugate points at all, while given any spacelike geodesic $\gamma$ there is no Jacobi field, with a spacelike covariant derivative along $\gamma$, that vanishes at two distinct points.

We say that a geodesic on $(\mathcal M,\textsl{g})$ is non-trapped if its maximal domain of definition is a bounded interval. Here the geodesic is allowed to intersect $\p \mathcal M$ in the interior points of its domain of definition.
Put differently, if $(\mathcal M,\textsl{g})$ is extended to a slightly larger manifold without boundary, then both ends of any maximally extended non-trapped geodesic must intersect the complement of $\mathcal M$. We assume
\begin{enumerate}
\item[(H4)] All null geodesics are non-trapped.
\end{enumerate}
This assumption is independent from (H1)--(H3) as can be seen for example by taking $\p \mathcal M$ to be the timelike hyperboloid $-t^2+|x|^2=1$ in the Minkowski spacetime $\R^{1+n}$. 

Finally, we assume 
\begin{enumerate}
\item[(H5)] All null geodesics have finite order of contact with the boundary.
\end{enumerate}
This technical condition is related to the geometric characterization of exact controllability of the wave equation on $\mathcal M$, see \cite{BG,BLR}.

We are now ready to state our main theorem. 

\begin{theorem}
\label{t0}
Let $(\mathcal M, \textsl{g})$ be a connected, smooth Lorentzian manifold with timelike boundary, and suppose that (H1)--(H5) hold. Let $V_1, V_2 \in C^\infty(\mathcal M)$.
If $\mathscr C(V_1) = \mathscr C(V_2)$ then $V_1 = V_2$.
\end{theorem}

In fact, we prove a slightly stronger result that is, roughly speaking, localized in the preimage $\{p \in \mathcal M : \tau(p) \in [-T,T] \}$
for large enough $T > 0$, see Theorem \ref{thm_main} below for the precise statement.

\subsection{On the curvature bound}
\label{curv_cond_sec}

Examples of non-constant curvature manifolds satisfying $R \le K$ can be constructed by considering Robertson--Walker type spaces
$(\mathcal M, \textsl{g})$ where 
$\mathcal M = \R \times M_0$ 
and $\textsl{g}$ has the warped product form
    \begin{align*}
-dt^2 + f(t) g_0(x), \quad (t,x) \in \R \times M_0.
    \end{align*}
Here the warping factor $f$ is smooth and strictly positive, and $(M_0, g_0)$ is a Riemannian manifold with boundary. 
If the sectional curvature of $(M_0, g_0)$ is bounded from above by a constant $C \in \R$, then, in view of \cite[Corollary 7.2]{AB08}, the manifold $(\mathcal M, \textsl{g})$ satisfies the bound $R\le K$ for any number $K$ in the interval 
$$\left[ \sup_{t\in \R} \frac{C+f'(t)^2}{f(t)^2}, \inf_{t \in \R} \frac{f''(t)}{f(t)}\right].$$

To see that there are manifolds satisfying (H1)--(H5)
such that their small perturbations satisfy the same conditions, we can simply take $f = 1$ and 
$(M_0, g_0)$ a Riemannian manifold with a strictly convex boundary and sectional curvature bounded from above by $C < 0$. This fact will be shown in detail in Section~\ref{sec_perturb}. Combining this example with Theorem~\ref{t0} we obtain the following corollary.
\begin{corollary}
\label{cor0}
Let $\mathcal M=\R\times M_0$, $g=-dt^2+g_0(x)$, where $(M_0,g_0)$ is a compact, simply connected Riemannian manifold with negative curvature and a smooth strictly convex boundary. Then given any metric $\widetilde{g}$ which is a sufficiently small perturbation of $g$ in $C^{\infty}_c(\mathcal M^{\textrm{int}})$, there holds:
$$\mathscr C(V_1) = \mathscr C(V_2) \implies V_1 = V_2.$$
Here, for $j=1,2$, $\mathscr C(V_j)$ is defined analogously to \eqref{C_V} corresponding to the metric $\widetilde{g}$ and potential $V_j \in C^{\infty}(\mathcal M)$.
\end{corollary}

\subsection{Previous literature}
\label{previous_lit}
Most of the previous results on the problem to find $V$ given $\mathscr C(V)$ are confined to the case of ultra-static cylinders, that is, $\mathcal M = \R \times M_0$
and $\textsl{g}$ is of the form (\ref{ortho_split})
with $c=1$ and $g_0$ independent from the time coordinate $t\in \R$. 
If also $V$ is assumed to be independent from $t$, then the problem is completely solved, without any further assumptions on $g_0$. 
The proof is based on the Boundary Control method that was introduced by Belishev in \cite{Bel87}. A geometric generalization to the case of ultrastatic cylinders is by Belishev and Kurylev in \cite{Bel92}. Conditional stability estimates were derived in \cite{AL04}. For later developments see the monograph \cite{KKL}, the review article \cite{Bel}, and the references in the recent paper \cite{KOP}. 

The Boundary Control method relies on the local unique continuation theorem by Tataru \cite{Tataru}. 
This result was subsequently generalized in \cite{RZ},
and its important precursors include the works of H\"ormander  \cite{Ho} and Robbiano \cite{Rob}.
Tataru's result allows for local unique continuation across any non-characteristic surface but requires that all the coefficients in the wave equation are real analytic with respect to the time variable. There is a broad class of counter examples due to Alinhac \cite{Al} showing that local unique continuation across non-pseudoconvex surfaces fails if smooth time-dependent coefficients are allowed. 
When both $\textsl{g}$ and $V$ are real analytic in the time coordinate, the result of Eskin \cite{E1,E2} solves the Lorentzian Calder\'{o}n problem using ideas from the Boundary Control method. 

Let us now discuss results that allow for an arbitrary $V \in C^\infty(\mathcal M)$.
Due to the gauge transformation (\ref{gauge}), there is no essential difference between the ultrastatic case and 
the case where $\textsl{g}$ is of the form (\ref{ortho_split})
with $g_0$ independent from the time coordinate
and with arbitrary $c$.
The first result with time-dependent $V$ is due to Stefanov in \cite{St} in the case that $(\mathcal M, \textsl{g})$ is the Minkowski spacetime. We also refer the reader to the work of Ramm and Sj\"ostrand \cite{RS} and Isakov \cite{I1} for similar results. 
The ultrastatic case was solved in \cite{FIKO}
under additional convexity conditions on the Riemannian factor $(M_0, g_0)$. The known results \cite{DKSU,DKLS} on the Riemannian Calder\'{o}n problem assume that $\textsl{g}$ is of the form $dt^2 + g_0(x)$. This can be viewed as the Riemannian analogue of the ultrastatic case.

The results \cite{St, RS, I1, FIKO} are based on a reduction of the Lorentzian Calder\'{o}n problem to the study of injectivity of the 
light ray transform $\mathscr LV$. 
Here $\mathscr LV(\gamma)$ is the integral of $V$ over $\gamma$, for inextendible null geodesics $\gamma$ on $\mathcal M$. 
As explained in \cite{SY} this reduction works for a very broad class of Lorentzian manifolds, however, inversion results for the light ray transform are scarce  outside the ultrastatic case. 
Stefanov inverted $\mathscr L$ in the case that $(\mathcal M, \textsl{g})$ is real analytic and satisfies a certain convexity condition \cite{Stefanov1}.
Recently $\mathscr L$ was inverted also in the case that $(\mathcal M, \textsl{g})$ is stationary, assuming additional convexity \cite{FIO}. While stationary metrics need not be the form (\ref{ortho_split}), they have time-independent coefficients. 

To summarize, in all the previous results $\textsl{g}$ has real analytic features. Most typically it is simply independent of the time variable, whereas in \cite{Stefanov1} it is real analytic, and in \cite{E1,E2} real analytic in the time variable. 
The main novelty of our result is that it requires no real analytic features. In fact, by Corollary~\ref{cor0}, the set of Lorentzian metrics satisfying (H1)--(H5) has a non-empty interior in the sense of arbitrary, smooth perturbations of the metric. We achieve this by introducing a novel controllability method, inspired by the Boundary Control method, that relies on a new unique continuation theorem in the exterior of the double null cone (Theorem~\ref{thm_uc}). A similar unique continuation theorem was previously known to hold in the Minkowski spacetime \cite{Alexakis,Shao}. 

We would like to emphasize that the curvature condition (H2) is used only in the proof of the unique continuation theorem. If unique continuation in the exterior of the double null cone is proven under weaker assumptions, then our controllability method will give stronger results on the Lorentzian Calder\'{o}n problem. We conjecture that the unique continuation theorem should hold under (H1) and (H3) only, and leave this as a direction for future research.   

Finally, we mention that the Lorentzian Calder\'{o}n problem has been
solved for non-linear wave equations in great geometric generality.
The version to determine a metric tensor g was solved, up to a
conformal factor, in \cite{KLU} for globally hyperbolic manifolds, with
data given on a small set in the interior rather than on  boundary. A
similar approach was used in \cite{FO} to determine a zeroth order
perturbation $V$, and the case with data on the boundary was considered
in \cite{HUZ}. However, all these works use non-linear terms in the wave
equation in an essential way, and their techniques do not generalize
for the linear wave equation considered in the present paper.

\subsection{Outline of the key ideas}

The key new ideas in this paper are outlined as follows. First, we introduce a novel {\em optimal} unique continuation result (Theorem~\ref{thm_uc}) for the wave equation in the exterior of the double null cone $\mathscr E_p$ from the point of view of the region where vanishing of the wave is \emph{assumed} and where vanishing is \emph{derived}. This optimality of the unique continuation result is absolutely essential in the proof of Theorem~\ref{t0}, in the same way that the optimal unique continuation result of Tataru was essential in the earlier works on wave equations with time-independent coefficients based on the Boundary Control method.  

Once Theorem~\ref{thm_uc} is proved, we use the result to obtain a novel controllability method from boundary, in the same spirit as that of Belishev, but with important differences. For instance, we do not posses here the typical boundary integral identities used in the Boundary Control method, namely the Blagovestchenskii identity \cite{Blago}. We use our unique continuation theorem together with an exact controllability result for the wave equation in rough Sobolev spaces to construct distributional solutions that, when restricted to a Cauchy surface, are supported at a point. We remark that the idea of using focusing waves has appeared in the literature, see \cite{Lassas} for a review. Finally, using these focused solutions together with special solutions (Gaussian beams), we are able to recover point values of the potential $V$ everywhere.   

\section{Strengthened main result}
\label{sec_tech_results}


Let $(\mathcal M,\textsl{g})$ be a connected, smooth Lorentzian manifold with signature $(-,+,\ldots,+)$ and a timelike boundary. Suppose that (H1) holds and let $T>0$. As we show in Appendix~\ref{appendix A}, the manifold admits an isometric embedding $\Phi:M\to \mathcal U\subset \mathcal M$ with $ \{q\in \mathcal M\,:\,\tau(q)\in [-T,T]\}\subset \mathcal U$, and 
\begin{equation}\label{cylinder}
M=[-T,T]\times M_0
\end{equation}
for some smooth, compact, connected manifold $M_0$ with smooth boundary and with the metric $g:=\Phi^*\textsl g$ satisfying
\begin{equation}\label{splitting} g(t,x)= c(t,x)\left(-dt^2+ g_0(t,x)\right),\quad \forall\, (t,x)\in M.\end{equation}
Here, $c$ is a smooth strictly positive function on $M$ and $g_0(t,\cdot)$ is a family of smooth Riemannian metrics on $M_0$ that depend smoothly on the variable $t\in [-T,T]$. As usual, $M_0=M_0^{\text{int}}\cup \p M_0$. 

In this section, we will formulate our unique continuation result and a stronger version of Theorem~\ref{t0} that are both stated on manifolds $(M,g)$ of the form \eqref{cylinder}--\eqref{splitting}. We note that the assumptions (H2)--(H3) and (H5) on $(\mathcal M,\textsl g)$ can also be formulated on $(M,g)$ in the natural way.

\subsection{Optimal unique continuation result}

We prove the following theorem in Section~\ref{uc_sec}.
\begin{theorem}\label{thm_uc}
Let $M=[-T,T]\times M_0$ be a Lorentzian manifold with a metric $g$ of the form \eqref{splitting} where $M_0$ is a smooth, compact, connected manifold with a smooth boundary. Suppose that $(M,g)$ satisfies (H2)--(H3) and (H5). Let $p \in M^{\text{int}}$ be such that $\xbar{\mathscr E_p}\subset (-T,T)\times M_0$, where $\mathscr E_p$ is defined by \eqref{exterior_null} with $\mathcal M=M$. Let $V \in C^{\infty}(M)$, $u \in H^{-s}(M)$ for some $s \geq 0$ and suppose that $(\Box+V) u=0$ on $M$ and also that the traces $u$ and $\p_\nu u$ both vanish on the set $\Sigma \cap  \mathscr E_p$, where $\Sigma=(-T,T)\times \p M_0$. Then, $u = 0$ on $\mathscr E_p$.
\end{theorem}

Note that since $\Sigma$ is a non-characteristic hypersurface for the wave operator, given any distributional solution to $(\Box+V)u=0$ on $M$, the traces of the distribution and all its derivatives on $\Sigma$ always exist. This property is called {\em partial hypoellipticity} (see for example, \cite[Proposition 7.6]{E3}). Similarly, the traces will be well-defined on $\{t\}\times M_0$ for any $t \in [-T,T]$. We remark also that the notation $H^{-s}(M)$ stands for the topological dual of $H^s_0(M)$ which itself is defined as the completion of smooth compactly supported functions in the interior of $M$, with respect to the standard $H^s(M)$-norm. 

Our proof of Theorem~\ref{thm_uc} is based on combining the principle of propagation of singularities for the wave equation (Lemma~\ref{lem_smoothness}), a layer stripping argument and a Carleman estimate with a degenerate weight, see equation \eqref{carleman_eq}. We remark that the unique continuation result here does goes beyond the observability estimates of Bardos, Lebeau and Rauch in \cite{BLR} and the unique continuation result of H\"ormander (for example \cite[Theorem 28.4.3]{Ho4}) in two key respects.
Firstly, as opposed to \cite{Ho4} where strict pseudoconvexity must be \emph{assumed} 
to derive  a uniqueness result, here a geometric condition on the curvature is identified which 
\emph{ensures} that the (Lorenzian, spacelike) distance function from a point  
satisfies a spacetime convexity property. The distance function is not strictly pseudoconvex, however. Secondly, as discussed before, contrary to both of these earlier results, our unique continuation result is optimal in terms of where vanishing of the wave is assumed and where it is derived.

We emphasize also that our unique continuation result does not follow from Tataru's unique continuation result \cite{Tataru} as the wave operator here has general time-dependent coefficients.

\subsection{Reduction from Cauchy data set to the Dirichlet-to-Neumann map on compact time intervals} We will show in this section that the Cauchy data set $\mathscr C(V)$ determines the Dirichlet-to-Neumann map for the wave equation on compact subsets $\{q\in \mathcal M\,:\,\tau(q)\in [-T,T]\}$ for all $T>0$. 

To make this precise, we begin by defining the Dirichlet-to-Neumann map $\Lambda_V$, associated to a manifold $M$ as in \eqref{cylinder} with a metric $g$ of the form \eqref{splitting} and a function $V \in C^{\infty}(M)$. Consider the wave equation 
\begin{equation}\label{pf0}
\begin{aligned}
\begin{cases}
\Box u+Vu=0\,\quad &\text{on $M$},
\\
u=f\,\quad &\text{on $\Sigma=(-T,T)\times \p M_0$,}\\
u(-T,x)=\p_{t}u(-T,x)=0 \,\quad &\text{on $M_0$,}
\end{cases}
    \end{aligned}
\end{equation}
where the wave operator, $\Box$, is defined in local coordinates $(t=x^0,\ldots,x^n)$ through the expression
$$ \Box u= -\sum_{j,k=0}^n \left|\det g \right|^{-\frac{1}{2}}\frac{\p}{\p x^j}\left(  \left|\det g \right|^{\frac{1}{2}}g^{jk}\frac{\p u}{\p x^k} \right).$$
Given each $f \in H^1_0(\Sigma)$, equation \eqref{pf0} admits a unique solution $u$ in the energy space
\begin{equation}
\label{natural_space}
C(-T,T;H^1(M_0))\cap C^1(-T,T;L^2(M_0)).
\end{equation}
Moreover, $\p_\nu u|_{\Sigma} \in L^2(\Sigma)$ where $\nu$ is the outward unit normal vector field on $\Sigma$, see for example \cite[Theorem 4.1]{LLT}. 

We define the Dirichlet-to-Neumann map, $\Lambda_{V}: H^1_0(\Sigma) \to L^2(\Sigma),$ by
\begin{equation}
\label{DNmap}
\Lambda_{V} f= \p_\nu u\big|_{\Sigma},
\end{equation}
where $u$ is the unique solution to \eqref{pf0} subject to the boundary value $f$ on $\Sigma$.  We remark that the definition of $\Lambda_{V}$ implicitly depends on the geometry $(M,g)$. We will however hide this dependence when no confusion arises.

Let us now consider a smooth, connected Lorentzian manifold $(\mathcal M,\textsl g)$ with a timelike boundary and suppose that it satisfies (H1)--(H5). We aim to show that $\mathscr C(V)$ determines the Dirichlet-to-Neumann map corresponding to certain subsets of $\mathcal M$. We start with a lemma.

\begin{lemma}
\label{lem_geo_tech}
Let $(\mathcal M,\text g)$ be a smooth connected Lorentzian manifold satisfying (H1) and (H4). Then, given any compact set $\mathcal K \subset \mathcal M$, there exists $T_{\mathcal K}>0$ such that 
$$\bigcup_{p\in \mathcal K}\xbar{\mathcal E_p} \subset \{q\in \mathcal M\,:\,\tau(q)\in (-T_{\mathcal K},T_{\mathcal K})\}.$$
\end{lemma}

\begin{proof}
We write $L^+ \mathcal M$ and $L^-\mathcal M$ for the bundle of future and past pointing null vectors on $T\mathcal M$ respectively, and write $L\mathcal M=L^+\mathcal M\cup L^-\mathcal M$. Given each $v=(q,\xi) \in L\mathcal M$ we define $\gamma_v:I\to \mathcal M$ to be the inextendible null geodesic with initial data $v$, that is to say,
\begin{equation}
\label{affine}
\nabla_{\dot{\gamma}_v(s)}\dot{\gamma}_v(s)=0,\quad \forall\, s \in I,
\end{equation}
subject to
$$ \gamma_v(0)=p \quad \text{and}\quad \dot{\gamma}_v(0)=\xi.$$
We write also $\gamma_{q,\xi}=\gamma_v$. Next, we define for each $v \in L^+\mathcal M$ the exit functions
\begin{align*}
 R^+(v)&=\sup \{s\in [0,\infty)\,:\, \gamma_v(s) \in M\}\quad \text{for each $v \in L^+\mathcal M$},\\
R^-(v)&=\inf \{s\in (-\infty,0]\,:\, \gamma_v(s) \in M\}\quad \text{for each $v \in L^+\mathcal M$}.
\end{align*}
Fixing an auxiliary Riemannian metric on $\mathcal M$, we write $S\mathcal M$ for the unit sphere bundle with respect to the auxiliary metric. It is straightforward to show that $R^+$ and $R^-$ are upper semi-continuous and lower semi-continuous respectively. We note also that $\tau\circ \gamma_{p,\xi}$ is a continuous and increasing function and that the composition of an upper (lower) semi-continuous functions with an increasing continuous function is upper (lower) semi-continuous.  Since $S\mathcal K\cap L^+\mathcal K$ is compact it follows that the map $(p,\xi) \mapsto \tau(\gamma_{p,\xi}(R^{\pm}(p,\xi))$ has a maximum (minimum) on $ L^+\mathcal K\cap S\mathcal K$ respectively.
\end{proof}

We have the following proposition.
\begin{proposition}
\label{prop_DNmap}
Let $(\mathcal M,\textsl{g})$ be a smooth, connected Lorentzian manifold with a timelike boundary and suppose that (H1)--(H5) hold. Let $V \in C^{\infty}(\mathcal M)$. Then, given any $T>0$, the Cauchy data set $\mathscr C(V)$ uniquely determines the  map $\Lambda_{V\circ \Phi}$ on $(M,g)$, corresponding to any fixed isometry $\Phi:M\to \mathcal U,$ 
with $\{q\in \mathcal M\,:\, \tau(q)\in [-T,T]\}\subset \mathcal U$, and with $M$ and $g=\Phi^*\textsl g$ of forms \eqref{cylinder} and \eqref{splitting} respectively.
\end{proposition}

\begin{proof}
Let
\begin{equation}
\label{f_space_1}
f \in C^{\infty}_c(\{q\in \p\mathcal M\,:\,\tau(q)\in (-T,T)\}). 
\end{equation}
We claim that there exists a unique element $(f,h)\in\mathscr C(V)$ with 
\begin{equation}
\label{h_space_1}
h\in C^{\infty}(\p\mathcal M)\quad \text{and}\quad \supp h \subset \{q\in \p\mathcal M\,:\,\tau(q)\in (-T,\infty)\}. 
\end{equation}
First, we show uniqueness. Suppose that $h_1, h_2 \in \mathscr C(V)$ satisfy \eqref{h_space_1}. For $j=1,2,$ let $u_j\in C^{\infty}(\mathcal M)$ solve $(\Box+V)u_j=0$ on $\mathcal M$ subject to $(u_j,\p_\nu u_j)|_{\p \mathcal M}=(f,h_j)$ and define $v=u_1-u_2$. Observe that $v|_{\p \mathcal M}=0$ and also that $\p_\nu v$ vanishes on $\{q\in \p\mathcal M\,:\,\tau(q)\in (-\infty,-T)\}$. 

By Lemma~\ref{lem_geo_tech}, there exists $T_1>0$ such that 
$$\bigcup_{q\in \tau^{-1}(-T)} \xbar{\mathscr E_q} \subset \{q\in \mathcal M\,:\, \tau(q) \in (-T_1,T_1)\}.$$
Let $p \in \mathcal M^{\text{int}}$ be chosen such that $\tau(p)<-T_1$ and observe that the set $\xbar{\mathscr E_p}$ is contained in $\{\tau<-T\}$. Applying Lemma~\ref{lem_geo_tech} again, we choose $\widetilde{T}>T$ sufficiently large so that $\xbar{\mathscr E_p} \subset \{|\tau|<\widetilde T\}$. Next, consider an isometry
$$\widetilde \Phi:[-\widetilde T,\widetilde T]\times \widetilde{M_0}\to \widetilde{\mathcal U},$$
with $\{|\tau|\leq \widetilde{T}\}\subset \widetilde{\mathcal U},$ and so that \eqref{splitting} holds. 

Writing $\widetilde{\Phi}^{-1}(p)=(T',x_0)$ and applying Theorem~\ref{thm_uc} on $[-\widetilde T,\widetilde T]\times \widetilde{M_0}$ for the function $v\circ \widetilde{\Phi}$, it follows that $$v(\widetilde{\Phi}(T',x))=\p_t v(\widetilde{\Phi}(T',x))=0 \quad \text{for all $x \in \widetilde{M_0}\setminus \{x_0\}$},$$
and consequently, by smoothness of $v$, $v(\widetilde{\Phi}(T',x))=\p_t v(\widetilde{\Phi}(T',x))=0$ for all $x \in \widetilde{M_0}$. Together with the fact that $(\Box +V)v=0$ on $\mathcal M$ and that $v|_{\p\mathcal M}=0$, it follows that $v$ must identically vanish everywhere in $\mathcal M$ and therefore $h_1=h_2$. 

To show existence, observe that given any $f$ as in \eqref{f_space_1}, there exists a unique smooth solution $u$ to the equation 
\begin{equation}\label{pf_infinite}
\begin{aligned}
\begin{cases}
\Box u+Vu=0\,\quad &\text{on $\mathcal M$},
\\
u=f\,\quad &\text{on $\p\mathcal M$,}\\
u=0 \,\quad &\text{on $\{\tau<-T\}$,}
\end{cases}
    \end{aligned}
\end{equation}
This is classical and follows for example from \cite[Chapter 7, Theorem 6]{Evans}. Setting $h=\p_\nu u|_{\p\mathcal M}$, it follows that $(f,h) \in \mathscr C(V)$ and that $h$ satisfies \eqref{h_space_1}.

Now, to finish the proof, we consider for each $f$ satisfying \eqref{f_space_1}, the unique $(f,h) \in \mathscr C(V)$ with $h$ satisfying \eqref{h_space_1}. Let $\Phi$ be as in the statement of the proposition. Note that $f\circ \Phi \in C^{\infty}_c(\Sigma)$ and that 
$$ \Lambda_{V\circ \Phi}(f\circ \Phi) =(h\circ\Phi)|_{\Sigma},$$
where $\Sigma=(-T,T)\times \p M_0$. The claim follows since $C^{\infty}_c(\Sigma)$ is dense in $H^1_0(\Sigma)$.
\end{proof}

\subsection{Local-in-time formulation of the main result}
We will prove in Section~\ref{thm_sec} the following slightly stronger version of Theorem~\ref{t0} that uniquely recovers a coefficient $V$ based on the knowledge of the Dirichlet-to-Neumann map $\Lambda_{V}$. 
\begin{theorem}
\label{thm_main}
Let $T>0$ and let $M=[-T,T]\times M_0$ be a Lorentzian manifold with a metric $g$ of the form \eqref{splitting} where $M_0$ is a smooth, compact, connected manifold with a smooth boundary. Suppose that $(M,g)$ satisfies (H2)--(H3) and (H5). Let $|T_0|<T$ and suppose that there is $p_0 \in (-T,T_0)\times M_0^{\text{int}}$ satisfying $\xbar{\mathscr E_{p_0}}\subset(-T,T_0)\times M_0$ where $\mathscr E_{p_0}$ is defined by \eqref{exterior_null} with $\mathcal M=M$. Let $V_1,V_2 \in C^{\infty}(M)$. Then, 
$$\Lambda_{V_1}f=\Lambda_{V_2}f \quad \forall\, f \in H^1_0(\Sigma) \implies V_1=V_2 \quad \text{on $\mathbb D$,}$$
where 
\begin{equation}\label{recovery_D} \mathbb D=\{p \in (T_0,T)\times M_0^{\textrm{int}} \,:\, \xbar{\mathscr E_p} \subset (T_0,T)\times M_0\}.\end{equation}
\end{theorem}

In the remainder of this section we show that the global Theorem~\ref{t0} follows immediately from the local version, Theorem~\ref{thm_main}. Indeed, to uniquely recover $V$ at an arbitrary point $p \in \mathcal M^{\text{int}}$ from the Cauchy data set $\mathscr C(V)$, we consider such $T_0\in \R$ that 
$$\xbar{\mathscr E_p}\subset \{q\in \mathcal M\,:\,\tau(q)\in (T_0,\infty)\}.$$ 
Note that this is possible thanks to Lemma~\ref{lem_geo_tech}. Next, we use the same lemma to choose $T_1>0$ such that
$$\bigcup_{q\in \tau^{-1}(T_0)} \xbar{\mathscr E_{q}}\subset\{q \in \mathcal M\,:\,\tau(q)\in (-T_1,T_1)\}.$$ 
We let $p_0\in \mathcal M^{\text{int}}$ with $\tau(p_0)< -T_1$ and observe that 
$$\xbar{\mathscr E_{p_0}}\subset\{q \in \mathcal M\,:\,\tau(q)\in (-\infty,T_0)\}.$$ 
Finally, we let $T>0$ be sufficiently large so that both the sets $\xbar{\mathscr E_p}$ and $\xbar{\mathscr E_{p_0}}$ are contained in $\{q \in \mathcal M\,:\,\tau(q)\in (-T,T)\}$.
Let $\Phi:M\to \mathcal U$ be an isometry where $\mathcal U$ contains $\{q\in \mathcal M\,:\, \tau(q)\in [-T,T]\},$ and with $M$ and $g=\Phi^*\textsl g$ satisfying \eqref{cylinder} and \eqref{splitting} respectively. By Proposition~\ref{prop_DNmap}, $\mathscr C(V)$ determines the map $\Lambda_{V\circ\Phi}$ on $(M,g)$. By invoking Theorem~\ref{thm_main} we deduce that $V(p)$ is determined by $\mathscr C(V)$. 

\subsection{Organization of the paper}

The remainder of this paper is organized as follows. In Section~\ref{sec_perturb} we give additional details on the perturbation example discussed in Section~\ref{curv_cond_sec} as well as other examples of geometries, with no underlying symmetries or real analyticity, that satisfy (H1)--(H5). In Section~\ref{prelim_sec}, we show that under the geometric hypotheses (H2)--(H3), there exists a smooth spacetime convex function in the exterior of the double null cone $\mathscr E_p$ emanating from a point $p$ in $M$. Section~\ref{uc_sec} contains the proof of the unique continuation theorem, namely Theorem~\ref{thm_uc}. In Section~\ref{ctrl_sec} we use the unique continuation theorem together with an observability estimate and a transposition argument to obtain exact controllability in rough Sobolev spaces, see Proposition~\ref{ctrl_prop}. Finally, in Section~\ref{thm_sec} we prove Theorem~\ref{thm_main} via a controllability method from the boundary. The controllability method allows us to construct waves that focus at a single point. Using focusing waves we recover point values of $V$ in $M$.

\section{Examples of manifolds satisfying (H1)--(H5)}
\label{sec_perturb}

The aim of this section is to illustrate that there are many Lorentzian manifolds satisfying geometric hypotheses (H1)--(H5) with no underlying symmetries or real analyticity. We begin by expanding on the perturbative example that was discussed in Section~\ref{curv_cond_sec} and show that Corollary~\ref{cor0} indeed follows from Theorem~\ref{t0}. Next, we discuss further examples of Lorentzian manifolds satisfying the curvature bound (H2) that is based on a local construction of metrics in a neighborhood of a fixed point in $\R^{1+n}$.
\subsection{Perturbations of ultrastatic manifolds with negative curvature}
Let $\mathcal M = \R \times M_0$ with the ultrastatic Lorentzian metric $g(t,x)=-dt^2+g_0(x)$,
where $(M_0, g_0)$ is a compact, simply connected Riemannian manifold of negative curvature with a smooth strictly convex boundary. 
We want to show that a small enough perturbation of the metric $g$ in $C^\infty_c(\mathcal M^\textrm{int})$ preserves (H1)--(H5). 

We will focus on (H2)--(H3) as showing the claim for (H1) and (H4)--(H5) is straightforward. We start with (H2) and write $\widetilde g$ for the perturbed metric. To simplify the notation, we also write $\epsilon>0$ for a generic small constant that roughly denotes the size of the perturbation. Recall that in view of \cite[Corollary 7.2]{AB08}, the unperturbed manifold $(M,g)$ satisfies
\begin{equation}\label{bound_curv_int} g(R(X,Y)Y,X)\leq K Q(X,Y),\end{equation}
for all $K \in [C,0]$ and $X,Y \in T_pM$ and $p \in \mathcal M$, where
$$ Q(X,Y)= g(X,X)g(Y,Y)-g(X,Y)^2.$$
We also note that given any $p \in \mathcal M$, any non-zero null vector $N \in T_p \mathcal M$ and any non-zero spacelike $Y \in T_p\mathcal M$ with $g(N,Y)=0$, there holds
\begin{equation}\label{bound_curv_null} g(R(N,Y)Y,N) <0.\end{equation}
Indeed, the latter inequality follows from the time-independence of the components of the metric $g$ and the fact that $(M_0,g_0)$ has a strictly negative curvature. To verify that (H2) is satisfied for small perturbations $\widetilde g$ of $g$ in the space $C^{\infty}_c(\mathcal M^{\textrm{int}})$ it suffices to prove the following stability lemma for the curvature bound. 
\begin{lemma}
\label{stability_curv_bound}
Let $(\mathcal M,g)$ satisfy the curvature bound $R \le K$ for all $K$ in some closed interval $[K_1,K_2]$. Assume also that the bound \eqref{bound_curv_null} holds. Then given any $K \in (K_1,K_2)$ and any compact set $U \subset \mathcal M$, there is a neighborhood $\mathscr G_U$ of $g$ in $C^\infty_c(U)$ such that for any $\widetilde g \in \mathscr G_U$ there holds $\widetilde R \le K$. Here $\widetilde R$ denotes the curvature tensor on $\mathcal M$ with respect to the perturbed metric $\widetilde g$.
\end{lemma}

\begin{proof}
We write $\widetilde Q$ analogously to $Q$ with $g$ replaced by $\widetilde g$ and write $\widetilde{\textrm{Sec}}$ for the sectional curvature with respect to $\widetilde g$. We also denote by $\epsilon$ the size of the small perturbation $\widetilde g$ with respect to $g$ in $C^{\infty}_c(U)$. Let $p\in \mathcal M$ be arbitrary and let $(x^0,\ldots,x^n)$ be a normal coordinate system centered at $p$ with $\frac{\p}{\p x^0}$ timelike. Given any $X,Y \in T_p\mathcal M$ we write 
$$X=\sum_{j=0}^nX^j \frac{\p}{\p x^j}\quad \text{and}\quad Y=\sum_{j=0}^nY^j \frac{\p}{\p x^j}.$$ 

As $\widetilde{\textrm{Sec}}(X,Y)$ only depends on $\Pi=\spn{\{X,Y\}}$, it suffices to show the claim for vectors $X$ and $Y$ satisfying 
\begin{equation}\label{vector_iden_1} X= c\frac{\p}{\p x^0} + Z\quad \text{and}\quad \widetilde{g}(Y,\frac{\p}{\p x^0})=\widetilde{g}(Z,\frac{\p}{\p x^0})=\widetilde{g}(Y,Z)=0,\end{equation}
with the scaling of the vectors normalized so that 
\begin{equation}\label{vector_iden_2}c^2+|Z|^2=|Y|^2=1,\end{equation}
where $|Z|^2=\sum_{j=1}^n |Z^j|^2$ and $|Y|^2=\sum_{j=1}^n |Y^j|^2$ denote Euclidean lengths of $Y$ and $Z$.
Observe that 
\begin{equation}\label{vector_iden_3} \widetilde{Q}(X,Y)=\widetilde{g}(X,X)\widetilde{g}(Y,Y)-\widetilde{g}(X,Y)^2=|Z|^2-c^2.\end{equation}

Given any $X$ and $Y$ as above, we first note that the claim follows easily if $|\widetilde Q(X,Y)|>\delta$ for $\delta>0$ small depending on the fixed choice of $K \in (K_1,K_2)$ and all $\epsilon$ sufficiently smaller than $\delta$. Indeed, away from degenerate planes (that is when $Q(X,Y)=0$) the sectional curvature changes by a small uniform amount depending on $\epsilon$. 

It remains to consider the case when the plane $\Pi=\spn{\{X,Y\}}$ is almost degenerate, that is to say $|\widetilde Q(X,Y)|<\delta$. Using \eqref{vector_iden_2}--\eqref{vector_iden_3} it follows that $|c|=\frac{\sqrt{2}}{2}+O(\delta)$, $|Z|=\frac{\sqrt{2}}{2}+O(\delta)$. Combining with \eqref{vector_iden_1}, we can write $X= N + \delta X'$
with $N$ a null vector satisfying $\sum_{j=0}^n |N^j|^2=1$ and  $\widetilde{g}(Y,X')=\widetilde{g}(Y,N)=0$ and $\sum_{j=0}^n|X'^j|^2<1$. As a consequence, given any $\epsilon$ sufficiently smaller than $\delta$ and in view of \eqref{bound_curv_null}, there holds
\begin{equation}\label{perturb_numerator}\widetilde{g}(\widetilde{R}(X,Y)Y,X)=g(R(N,Y)Y,N)+O(\delta)\in (C_1,C_2),\end{equation}
for some constants $C_1<C_2<0$ independent of $\delta$. 

We have two cases depending on the sign of $\widetilde Q(X,Y)$. First, let us consider the case that the plane is timelike and almost degenerate, that is to say 
\begin{equation}\label{degenerate_Q}\widetilde Q(X,Y) \in (-\delta,0).\end{equation} 
Then, applying \eqref{perturb_numerator} we conclude that
$$ \widetilde{\textrm{Sec}}(X,Y)=\frac{\widetilde{g}(\widetilde{R}(X,Y)Y,X)}{\widetilde{Q}(X,Y)}>-\frac{C_1}{2\delta}>K.$$
Next, we consider the remaining case, that is when $Q(X,Y) \in (0,\delta)$. Now, it follows from \eqref{perturb_numerator} that
$$ \widetilde{\textrm{Sec}}(X,Y)=\frac{\widetilde{g}(\widetilde{R}(X,Y)Y,X)}{\widetilde{Q}(X,Y)}<\frac{C_2}{2\delta}<K.$$
\end{proof}
\begin{remark}
We note that the proof above works also for small $C^2(\mathcal M)$ perturbations of the metric $g$ with compact support in $U$. For the sake of simplicity, we will not keep track of the precise space $g \in C^{k}(\mathcal M)$ for which Theorem~\ref{t0} holds and instead simply work with metrics in $C^{\infty}(\mathcal M)$.  
\end{remark}
Next, we show that (H3) also holds for small enough perturbations $\widetilde g$ of the ultrastatic metric $g$ above.  To fix the notation, we write for each $p,q\in M$, $p \le q$ if there is a causal piecewise smooth path on $M$ from $p$  to $q$ or $p=q$. Also, we denote  $p\ll q$
if there is a future pointing piecewise smooth timelike path on $M$ from $p$ to $q$. Using these relations, the causal future and past of a point $p\in M$ can be defined by
\begin{equation}\label{causal_rel}J^+(p) = \{x \in M: p \le x\}\quad \text{and} \quad J^-(p) = \{x \in M : x \le p\}.\end{equation}
The chronological future and past of $p$ is defined analogously with the causal relation replaced by the chronological relation,
$$I^+(p)=\{x \in M: p \ll x\}\quad \text{and}\quad I^-(q) = \{x \in M : x \ll q\}.$$ 

Note that for the unperturbed manifold $(M,g)$, the exponential map $\exp_p$ is a diffeomorphism onto $\mathcal M$ for all $p \in \mathcal M$. Moreover, for all $T > 0$,
$\exp_p$ on $M := [-T,T] \times M_0$ 
is also a diffeomorphism onto $M$ for all $p \in M$. A small enough perturbation of the metric in $C^\infty_c(M^\textrm{int})$
preserves the latter fact for any fixed $T$. 
It also preserves the fact that all null geodesics hit $\p \mathcal M$ transversally.
\HOX{It might be more natural to say that $T \Sigma$ contains no null vectors.}

\HOX{Check this, should be easy since we will never see the perturbation as we work far away from it.}
Choosing large enough $T > 0$ in comparison to the size of the support of the perturbation, the perturbed manifold satisfies the following: there are no cut points on null geodesics not contained in $M$ and, for all $p \in \mathcal M$ such that $\mathscr E_p$ is not contained in $M$, $\exp_p$ is a diffeomorphism from the spacelike vectors in its maximal domain of definition onto $\mathscr E_p$.

To complete the perturbation argument, it suffices to show the following two lemmas. 
\begin{lemma}\label{lem_no_cutpts}
Suppose that $\exp_p$ on $(M,\widetilde g)$ 
is a diffeomorphism onto $M$ for all $p \in M$,
and that all null geodesics are transversal to $\Sigma$.
Let $\gamma$ be an inextendible null geodesic contained in $M$. 
Then for any points $p< q$ on $\gamma$,
the only causal path in $\mathcal M$ from $p$ to $q$ is along $\gamma$.
\end{lemma}

\begin{lemma}\label{lem_normalcoords}
Suppose that $\exp_p$ on $(M,\widetilde g)$ 
is a diffeomorphism onto $M$ for all $p \in M$,
and that all null geodesics are transversal to $\Sigma$.
Let $p \in M$ satisfy $\mathscr E_p \subset M$. 
Then the exponential map $\exp_p$
is a diffeomorphism from the spacelike vectors onto $\mathscr E_p$.
\end{lemma}

We emphasize that the causality relations and the set $\mathscr E_p$ in the preceding lemmas are defined with respect to the metric $\widetilde g$. The following lemma is a variant of \cite[Proposition 4.5.1]{HE} and we omit its proof.

\begin{lemma}\label{lem_simplecausality}
Let $p \in M$ and suppose that $\exp_p$ is a diffeomorphism onto $M$,
and that all null geodesics are transversal to $\Sigma$.
Then $I^+(p)$ is the image of future pointing timelike vectors under $\exp_p$,
and $J^+(p)$ is the image of future pointing causal vectors under $\exp_p$. The analogous statements hold for $I^-(p)$ and $J^-(p)$.
\end{lemma}

Note that Lemma \ref{lem_simplecausality} implies Lemma \ref{lem_normalcoords}. Indeed, $\mathscr E_p$ is the complement of $J^+(p) \cup J^-(p)$ and so it is the image of spacelike vectors. 

The following lemma is a variant of \cite[Proposition 10.46]{Oneill}. 

\begin{lemma}\label{lem_shortcut}
Suppose that all null geodesics are transversal to $\Sigma$.
Let $p, q \in M$.
If there is a causal path in $M$ from $p$ to $q$ that is not a null pregeodesic, then $p \ll q$.
\end{lemma}
\begin{proof}
Let $\gamma : [a,b] \to M$ be a causal from $p$ to $q$ that is not a null pregeodesic. Note that 
if $\gamma(s) \in \Sigma$ for $s$ in an open interval $I \subset [a,b]$ then $\dot \gamma(s)$
is timelike for $s \in I$ since $T \Sigma$ contains no null vectors.
Let us consider a maximal interval $[s_1, s_2] \subset [a,b]$
such that $\gamma(s) \in M^\textrm{int}$ for $s \in (s_1, s_2)$.
We will show first that there is timelike path joining $\gamma(s_1)$ and $\gamma(s_2)$ or $\gamma$ is a null pregeodesic on $[s_1, s_2]$. 

If $\dot\gamma(s_1)$ is timelike then there is $\tilde s_1 > s_1$ such that $\gamma$ is timelike on $[s_1, \tilde s_1]$ and $\gamma(\tilde s_1) \in M^\textrm{int}$.
In this case it is enough to show that there is a timelike path joining $\gamma(\tilde s_1)$ and $\gamma(s_2)$, and we replace $s_1$ by $\tilde s_1$.
We shorten the interval $[s_1, s_2]$ analogously if $\dot\gamma(s_2)$ is timelike.
If $\dot\gamma(s_j)$ is lightlike and $\gamma(s_j) \in \p M$ for $j=1$ or $j=2$,
then $\dot\gamma(s_j)$ is not tangential to $\Sigma$ since $T \Sigma$ contains no null vectors.
Thus we may assume without loss of generality that $\alpha := \gamma|_{[s_1,s_2]}$ is transversal to $\Sigma$.
If $\alpha$ is not a null pregeodesic, then the proof of [O'Neill, Prop. 10.46] applies to $\alpha$, and shows that $\alpha$ can be perturbed 
so that it becomes a timelike path from $\gamma(s_1)$ to $\gamma(s_2)$.

It remains to consider the case that $\alpha$
is a null geodesic and $\alpha(s_j) \in \Sigma$ for both $j=1$ and $j=2$. 
Then $s_1 > a$ or $s_2 < b$.
We consider only the former case, the latter being similar.
To simplify the notation we reparametrize $\gamma$ so that $s_1 = 0$.
We show that for small $\epsilon > 0$ the path
$\gamma|_{[-\epsilon, \epsilon]}$ can be perturbed so that it becomes timelike while keeping the end points $\gamma(-\epsilon)$ and $\gamma(\epsilon)$ fixed. This concludes the proof as it allows us to reduce to the case that $\alpha$
is not a null pregeodesic. 

To make the discussion explicit, let us consider boundary normal coordinates $(y,r) \in \Sigma \times [0,r_0)$ such that 
    \begin{align*}
\widetilde g(y,r) = \begin{pmatrix}
\widetilde h(y,r) & 0 \\
0 & 1
\end{pmatrix}.
    \end{align*}
Let $w$ be the inward pointing unit normal vector to $\Sigma$, that is, $w = (0,1)$ in the above coordinates. We write $W$ for the parallel translation of $w$ along $\gamma$ and set $V = fW$ where $f : [-\epsilon, \epsilon] \to [0,1]$ is a smooth function for some $\epsilon > 0$. We choose $f$ so that its derivative satisfies $f'(s) > 0$ for $s < 0$ and $f'(s) < 0$ for $s > 0$ and that $f(\pm \epsilon ) = 0$.
The covariant derivative of $V$ satisfies $V' = f' W$ since $W$ was obtained as a parallel translation.
We show that a small variation of $\gamma|_{[-\epsilon, \epsilon]}$, with the variation field $V$, is timelike. Note that this variation is well-defined in $M$ for small $\epsilon > 0$ since $f \geq 0$ and, writing $W = (W^0, \dots, W^n)$, there holds $W^n \geq 0$. Moreover, the variation keeps $\gamma(\pm \epsilon)$ fixed since $f(\pm \epsilon) = 0$.

Recall the notation $\alpha(s) = \gamma(s)$ for $s \geq 0$ near the origin, and write $\beta$ for the restriction of $\gamma$ in $s \leq 0$.
We have $\dot \alpha^n(0) > 0$ since $\alpha$ is transversal to $\Sigma$. Hence $\widetilde g({w, \dot \alpha(0)}) > 0$
and also $\widetilde g({W, \dot \alpha}) > 0$ on $[0, \epsilon]$ for small $\epsilon > 0$. Furthermore,  $\dot \beta^n(0) \leq 0$
with the equality possible only if $\dot \beta(0)$ is timelike. 

Let us consider first the case that $\dot \beta^n(0) < 0$.
Then $\widetilde g({W, \dot \beta}) < 0$ on $[-\epsilon, 0]$ for small $\epsilon > 0$. Taking into account the sign of $f'$, we see that $\widetilde g({V', \dot \gamma}) < 0$ for $s \in [-\epsilon, \epsilon]$. Now \cite[Lemma 10.45]{Oneill} implies that a small variation of $\gamma|_{[-\epsilon, \epsilon]}$, with the variation field $V$, is timelike as required. 

We turn to the case that $\dot \beta^n(0) = 0$ and $\dot \beta(0)$ is timelike. 
Then $\dot \beta$ is timelike on $[-\epsilon, 0]$ for small $\epsilon > 0$. On one hand, \cite[Lemma 10.45]{Oneill} implies that a small variation of $\alpha|_{[0, \epsilon]}$, with the variation field $V$, is timelike. On the other hand, any  small variation of $\beta = \gamma|_{[-\epsilon, 0]}$ is timelike. 
\end{proof}

\begin{proof}[Proof of Lemma \ref{lem_no_cutpts}]
To get a contradiction, suppose that there is a future pointing causal path $\beta$ from $p$ to $q$ for some $p< q$ on $\gamma$,
and that $\beta$ is not a reparametrization of $\gamma$.
Note that $\beta$ can not leave $M$ since the temporal function (that is, the coordinate on the factor $\R$ of $\mathcal M$)
is increasing on $\beta$.
The case that $\beta$ is a null pregeodesic is a contradiction with $\exp_p$ being an injection. 
Thus $\beta$ is not a null pregeodesic, and Lemma \ref{lem_shortcut} implies that $p \ll q$.
But now Lemma \ref{lem_simplecausality} implies that there is a timelike geodesic from $p$ to $q$, a contradiction with $\exp_p$ being an injection.
\end{proof}

\subsection{Further examples of manifolds satisfying $R\le K$}
We have already seen that perturbations of warped product spaces contain many examples of geometries satisfying (H1)--(H5). In this section we discuss a local method of constructing Lorentzian metrics in small neighborhoods of a fixed point in $\R^{1+n}$ that satisfy the curvature bound $R\le K$. This will further illustrate the richness of manifolds satisfying (H1)--(H5). To this end, let us begin by fixing the origin in $\R^{1+n}$ and defining for each index $i,j=0,\ldots,n$, the tensor $g_{ij}=g(\p_i,\p_j)$ by the expression
\begin{equation}\label{metric_example} g_{ij}=\eta_{ij} -\frac{1}{3}R_{iklj}(0)x^kx^l + O(|x|^3),\end{equation}
where $\eta_{ij}$ is the Minkowski metric and $R_{iklj}(0)$ is a rank four tensor at the origin that we will explicitly construct next. Our aim is to construct the tensor $R_{ijkl}(0)$ in such a way that the estimate $R\le K$ holds in a small neighborhood of the origin in $\R^{1+n}$. In fact, this can be done in many ways but we just show one such example. We remark that the choice of the symbol $R_{ijkl}(0)$ here is deliberate as it is well-known that given an metric $g$ of the form \eqref{metric_example}, the curvature tensor at the origin (evaluated in a normal coordinate system) will be equal to $R_{ijkl}(0)$. 

Writing $\alpha$, $\beta$, $\delta$ and $\gamma$ for indices running from $1$ to $n$, we define 
$$R_{0\alpha\beta0}(0)=R_{\alpha00\beta}(0)=-R_{\alpha0\beta0}(0)=-R_{0\alpha0\beta}(0)=C_{\alpha\beta},$$ 
where $C_{\alpha\beta}$ is an arbitrary symmetric matrix whose eigenvalues are bounded from above by zero. Next, we set 
$$R_{0\alpha\beta\delta}(0)=R_{\alpha0\beta\delta}(0)=R_{\alpha\beta0\delta}(0)=R_{\alpha\beta\delta0}(0)=0$$ 
and 
$$R_{00\alpha\beta}(0)=R_{\alpha\beta00}(0)=0.$$
Finally, we set
$$ R_{\alpha\beta\delta\gamma}(0)=\frac{\kappa}{2}(I\KN I)_{\alpha\beta\delta\gamma},$$
where $\wedge\!\!\!\!\!\!\!\!\!\!\;\;\bigcirc$ is the Kulkarni-Nomizu product, $\kappa>0$ is arbitrary and $I$ is the identity matrix. By definition of the Kulkarni-Nomizu product we obtain that
$$R_{\alpha\beta\delta\gamma}(0)X^{\alpha}Y^{\beta}Y^{\delta}X^{\gamma}=-\kappa \left((\sum_{j=1}^n |X^j|^2)(\sum_{j=1}^n|Y^j|^2)-(\sum_{j=1}^nX^jY^j)^2\right),$$
where we are using the Einstein summation convention with $\alpha$, $\beta$, $\delta$, $\gamma$ running between $1$ and $n$.

Note that to show the curvature bound at the origin, it suffices to show it for $X,Y \in T_0\R^{1+n}$ with $Y^0=0$. Now,
\begin{align*}
g(R(X,Y)Y,X)&=R_{ijkl}(0)X^iY^jY^kX^l\\
&= C_{\alpha\beta}|X^0|^2Y^{\alpha}Y^{\beta}-\kappa \left((\sum_{j=1}^n |X^j|^2)g(Y,Y)-g(X,Y)^2\right)\\
&\leq K\, Q(X,Y),
\end{align*}
for all $K \in [-\kappa,0]$, where in the last step we use the fact that $\kappa>0$ and $C_{\alpha\beta}$ is a strictly negative tensor. It can also be verified that given any non-zero null vector $N$ at the origin and any non-zero spacelike vector $Y$ orthogonal to $N$, the bound \eqref{bound_curv_null} holds at the origin. Thus, by stability (analogously to the proof of Lemma~\ref{stability_curv_bound}) and given each $K \in (-\kappa,0)$ the curvature bound $R\le K$ must hold in a sufficiently small neighborhood of the origin. 

\section{A distance function and its Hessian}
\label{prelim_sec}
\subsection{Notations}
\label{prelim_sec_not}

We start with fixing some notations for the remainder of the paper. We let $T>0$ and consider $M$ and $g$ as in \eqref{cylinder} and \eqref{splitting} respectively. We assume that the factor $M_0$ in \eqref{cylinder} is a connected, compact, smooth Lorentzian manifold with a smooth boundary. We also write $\Sigma=(-T,T)\times \p M_0$ for the timelike boundary of $M$.

Given each $q \in M$, we use the notation 
$$\langle X,Y\rangle=g(X,Y),\quad \forall \,X,Y \in T_qM,$$ 
and let $\nabla$ and $\div$ denote the gradient and divergence operator on $M$ respectively. We also define the Hessian of a function $\psi:M\to \R$ through
$$ \Hess \psi(X,Y) = \langle D_X \nabla \psi,Y\rangle \quad \forall\, X,Y \in T_qM,$$
where $D_XY$ is the Levi-Civita connection on $M$. 

For the functional spaces, we recall that given any smooth submanifold $U\subset M$ with smooth boundary, we use the notation $H^s(U)$, $s\geq 0$ to stand for the standard Sobolev spaces on $U$ and recall that $H^s_0(U)$ denotes the completion of smooth compactly supported functions in the interior of $U$, with respect to the norm $\|\cdot\|_{H^s(U)}$. We also recall, for each $s\geq 0$, the space $H^{-s}(U)$ that is defined as the topological dual of $H^s_0(U)$. Finally, we write $\mathcal E'(U)$ for the set of compactly supported distributions in the interior of $U$. 

Given two smooth compactly supported functions $v,w$ in the interior of $U$, we write
$$ (v,w)_{U} = \int_{U} v\,w\,dV_{g},$$
where $dV_{g}$ is the induced volume form on the submanifold $U$. Note that the latter pairing $(\cdot,\cdot)_U$ admits a continuous extension to $H^{-s}(U)\times H^s_0(U)$ for each $s \geq 0$. When no confusion is likely to arise, we will use $(\cdot,\cdot)_U$ in this extended sense.

\subsection{A distance function in the exterior of the double null cone}

Let $p \in M$ and let $e_0,\dots,e_n$ be an orthonormal basis of $T_p M$
in the sense of \cite[Lemma 24, p. 50]{Oneill}, that is,
for distinct $j,k = 0,\dots,n$,
    \begin{align*}
\pair{e_j, e_k} = 0, \quad
\pair{e_j, e_j} = \epsilon_j,
    \end{align*}
where $\epsilon_0 = -1$ and $\epsilon_j=1$ for $j=1,\ldots,n$.
We consider the function 
    \begin{align}\label{def_r}
r_p(y^0,\ldots,y^n) = \left( \epsilon_0(y^0)^2 + \dots + \epsilon_n (y^n)^2 \right)^{1/2}
    \end{align}
in the corresponding normal coordinates $y^0,\dots,y^n$ centered at the point $p$. Note that under condition (H3), the normal coordinate system is well defined on the set $\mathscr E_p$ defined by \eqref{exterior_null}. 

We will show that the function $r_p$  is a distance function on $\mathscr E_p$, that is to say $\pair{\nabla r_p, \nabla r_p} = 1$. To see this, consider the local hyperquadric 
    \begin{align*}
Q = \{\omega \in T_p M : \pair{\omega,\omega} = 1\}.
    \end{align*}
In the region $\{r_p > 0\}$, we consider the polar coordinates $x = r_p\omega$ with $r_p>0$ and $\omega \in Q$.
The Gauss lemma, see e.g. \cite[Lemma 1, p. 127]{Oneill}, implies that in these coordinates, with $r_p$ having the index $n$,
the metric tensor has the form
    \begin{align*}
\begin{pmatrix}
h(r_p\omega) & 0
\\
0 & 1
\end{pmatrix}.
    \end{align*}
It follows that $\pair{\nabla r_p, \nabla r_p} = 1$. We also have the following well-known lemma on distance functions in general. We give a short proof for the convenience of the reader.

\begin{lemma}
\label{lem_r}
Let $r$ be a distance function and denote its gradient by $\p_{r}$.
Then $\nabla_{\p_{r}} \p_{r} = 0$.
\end{lemma}
\begin{proof}
For a vector field $Z$,
    \begin{align*}
0 = Z1 = Z \pair{\p_{r}, \p_{r}} = 2 \pair{\nabla_Z \p_{r}, \p_{r}}
    \end{align*}
and hence
    \begin{align*}
0 = \Hess r(Z,\p_{r}) = \Hess r(\p_{r},Z) = \pair{\nabla_{\p_{r}} \p_{r}, Z}.
    \end{align*}
\end{proof}

\subsection{Spacetime convex functions}

In this section, we recall one of the results of the work of Alexander and Bishop in \cite{AB08} that relates the curvature bound in hypothesis (H2) with the existence of a spacetime convex function on the manifold. This function will subsequently be used in the proof of Theorem~\ref{thm_uc}. The notion of spacetime convex functions was first introduced by Gibbons and Ishibashi in \cite{GI} purely for geometrical pursuits and not related to any unique continuation results. A function $\psi:M\to\R$ is spacetime convex if it satisfies, at each point $p \in M$, the bound
$$ \Hess \psi(X,X) \geq \lambda\, \langle X,X\rangle\quad \text{for some $\lambda \in \R$ and all $X \in T_pM$}.$$ 

We remark that spacetime convexity differs from H\"ormander's strict pseudoconvexity (see \cite[Chapter 28]{Ho4}), which requires at each point $p \in M$ the bound
$$ \Hess \psi(X,X) >0,\quad \text{if $\langle X,X\rangle=\langle X,\nabla \psi\rangle=0$.}$$
Due to the strictness of the latter inequality, strict pseudoconvexity does not follow from spacetime convexity. Strict pseudoconvexity was used by H\"ormander, see \cite[Theorem 28.4.3]{Ho4}, to prove unique continuation for solutions to the wave equation locally near a point across level sets of $\psi$. We show in Section~\ref{uc_sec} that one can also derive unique continuation results using the alternative spacetime convexity.

To construct our spacetime convex function, we start with a definition.

\begin{definition}
Under the hypotheses (H2)--(H3) for $(M,g)$, we define for each $p \in M$, a function on $\mathscr E_p$ through
\begin{equation}\label{def_f}
\psi_{K,p}(q)=\begin{aligned}
\begin{cases}
\sqrt{|K|}r_p(q)\cot(\sqrt{|K|}r_p(q)), \quad &\text{if $K>0$.}
\\
1,\quad &\text{if $K=0$},
\\
\sqrt{|K|}r_p(q)\coth(\sqrt{|K|}r_p(q)), \quad &\text{if $K<0$.}
\end{cases}
    \end{aligned}
\end{equation}
\end{definition} 

We note that since $r_p(q)$ is positive on $\mathscr E_p$, it follows easily that when $K\leq 0$, the function $\psi_{K,p}$ is positive on $\mathscr E_p$. In the case $K>0$, the assumption $\textrm{Diam}(M)<\frac{\pi}{2 \sqrt{K}}$ in (H2) implies that $\cot(\sqrt{|K|}r_p(q))>0$ on the set $\mathscr E_p$.

\begin{lemma}
\label{spacetime_conv_lem}
Let $(M,g)$ satisfy (H2)--(H3) for some $K\in \R$. Let $p \in M$ and consider the function $\psi_{K,p}$ defined by \eqref{def_f}. Then,
$$\Hess r_p(X,X) \geq \frac{\psi_{K,p}(q)}{r_{p}(q)}\,\left(\langle X,X\rangle-\langle X,\nabla r_p\rangle^2\right),$$
where the Hessian is computed at a point $q \in \mathscr E_p$ and $X \in T_q \mathscr E_p$. 
\end{lemma}

\begin{proof}
We start by defining the function $\widetilde{\psi}_{K,p}:\mathscr E_p \to \R$ by the expression
\[
\widetilde{\psi}_{K,p}(q)=\begin{aligned}
\begin{cases}
\frac{1}{|K|}(1-\cos(\sqrt{|K|}r_p(q))), \quad &\text{if $K>0$.}
\\
\frac{1}{2}{r_p(q)^2},\quad &\text{if $K=0$},
\\
\frac{1}{|K|}(-1+\cosh(\sqrt{|K|}r_p(q))), \quad &\text{if $K<0$.}
\end{cases}
    \end{aligned}
\]
Under the hypotheses (H2)--(H3) and by applying \cite[Corollary 4.5--4.6]{AB08} (see also \cite[Theorem 4.6]{AK}), it follows that
$$\Hess{\widetilde{\psi}_{K,p}}(X,X) \geq (1-K\widetilde{\psi}_{K,p}(q))\langle X,X\rangle,\quad \forall\, X\in T_q\mathscr E_p.$$ 
The claim follows by rewriting the above inequality in terms of the distance function.
\end{proof}

\section{Unique continuation in the exterior of the double null cone}
\label{uc_sec}
\subsection{Smoothness away from the double null cone}

We start with a purely geometric lemma.

\begin{lemma}
\label{lem_cutlocus}
Let $M=[-T,T]\times M_0$ and let $g$ be a smooth Lorentzian metric of the form \eqref{splitting} where $M_0$ is a compact smooth manifold with smooth boundary. Suppose also that $M$ satisfies (H3). Let $p \in M^{\textrm{int}}$, $q \in \mathscr E_p$ and let $\gamma:I\to M$ to be an inextendible null geodesic passing through $q$. Then $\gamma$ intersects $J^+(p)\cup J^-(p)$ at most once.
\end{lemma}

\begin{proof}
Without loss of generality, we assume that $\gamma$ is future pointing. To get a contradiction, we assume that $\gamma$ intersect $J^+(p) \cup J^-(p)$
at two points $q_1 < q$ and $q_2 > q$.

In the case that $p \le q_1$ and $p \le q_2$
there is a causal path from $p$ to $q$ via $q_1$,
a contradiction with $q \in \mathscr E_p$.
The case that $p \ge q_1$ and $p \ge q_2$ leads to an analogous contraction. 
Finally, in case that $p \ge q_1$ and $p \le q_2$ there are two causal paths from $q_1$ to $q_2$,
one via $q$ along $\gamma$ and the other via $p$, a contradiction with (H3).
\end{proof}
The above lemma can be used together with the principle of propagation of singularities for solutions to the wave equation to deduce that solutions to the wave equation with vanishing Cauchy data in $\mathscr E_p \cap \Sigma$ must be smooth in $\mathscr E_p$. More precisely,
\begin{lemma}\label{lem_smoothness}
Let the Lorentzian manifold $(M,g)$ and the distribution $u$ be as in the hypothesis of Theorem~\ref{thm_uc}. Then,
$$u\in C^{\infty}(\mathscr E_p).$$
\end{lemma}

\begin{proof}
Let us consider a smooth extension of the manifold $M_0$ into a slightly larger manifold $\widetilde M_0$. We write $\widetilde M=[-T,T]\times \widetilde{M_0}$ and extend the metric $g$ smoothly to $\widetilde M$. Finally, we extend the function $u$ to all of $\widetilde M$ by setting it to zero on $\widetilde M\setminus M$ and write $\widetilde u$ for this extended function. 

Note that the traces $u$ and $\p_\nu u$ both vanish on the timelike hypersurface $\mathscr E_p\cap \Sigma$. Thus, $\Box \widetilde u = f$ on $\widetilde M$ where $f$ is a distribution with $\supp f \subset M \setminus \mathscr E_p$. Let $q \in \mathscr E_p$, $\xi \in L_qM$ and write $v=(q,\xi)$. We write $\widetilde{\gamma}_v$ for the inextendible null geodesic in $\widetilde M$ with initial data $v$. In view of Lemma~\ref{lem_cutlocus}, it follows that there is a point $\tilde{q}\in \widetilde M\setminus M$ on $\widetilde{\gamma}_v$ and a neighborhood $U\subset \mathscr E_p\cup (\widetilde M\setminus M)$ of the segment of $\widetilde{\gamma}_v$ between $q$ and $\tilde q$. Note that the cases $\tilde q<q$ and $q<\tilde q$ are both possible. As $\widetilde u=0$ near $\tilde q$ we can apply \cite[Theorem 26.1.4]{Ho4} to conclude that the wave front set of $u$ at $q$ does not contain $\xi$. As the null vector $\xi$ was arbitrary, it follows that the wave front set of $u$ at $q$ can not contain any null vectors. The claim follows thanks to the microlocal ellipticity of the wave equation in the non-null directions.
\end{proof}

\subsection{The conjugated wave operator}

In this section, we show the following pointwise identity. This lemma will in fact be true over general semi-Riemannian manifolds.
\begin{lemma}\label{lem_basic_carleman}
Let $U \subseteq M$ and let $v, \sigma \in C^2(\xbar{U})$ and $\ell \in C^3(\xbar{U})$.
Then
    \begin{align*}
(e^\ell \Box (e^{-\ell} v))^2 / 2
= S^2/2 + 2 \pair{\nabla \ell, \nabla v}^2 + P + R + \div B
    \end{align*}
where the squared part is
    \begin{align*}
S &= \Box v + q v, \quad q = -\pair{\nabla \ell, \nabla \ell} - \Box \ell - \sigma,
    \end{align*}
writing $\tilde \sigma = \sigma + \Box \ell$, the leading part is
    \begin{align*}
P = \tilde \sigma \pair{\nabla v, \nabla v} + 2 \Hess \ell(\nabla v, \nabla v)+ \left(-\tilde \sigma \pair{\nabla \ell, \nabla \ell} + 2 \Hess \ell(\nabla \ell, \nabla \ell) \right)v^2,
    \end{align*}
and the remainder and divergence parts are
    \begin{align*}
R &= \left(\div ((\Box \ell)\; \nabla \ell) - \sigma \tilde \sigma + \frac12 \sigma^2 + \frac12 \Box \sigma\right) v^2,
\end{align*}
and
\begin{multline*}
B= -(2\pair{\nabla \ell, \nabla v} + \sigma v) \nabla v+ \frac12 v^2 \nabla \sigma\\
+(\pair{\nabla v, \nabla v} - (\Box \ell + \pair{\nabla \ell, \nabla \ell}) v^2) \nabla \ell.
    \end{multline*}
\end{lemma}
\begin{proof}
We have
    \begin{align*}
e^\ell \Box (e^{-\ell} v)=
\Box v - \Box \ell\; v + 2 \pair{\nabla \ell, \nabla v}
- \pair{\nabla \ell, \nabla \ell} v
= S + A,
    \end{align*}
where $A = 2\pair{\nabla \ell, \nabla v} + \sigma v$. 
Thus,
    \begin{align}\label{conj_box_square}
(e^\ell  \Box (e^{-\ell} v))^2/2 = S^2/2 + A^2/2 + SA.
    \end{align}
Moreover,
    \begin{align}\label{A_squared}
A^2/2 = 2\pair{\nabla \ell, \nabla v}^2 + \sigma^2 v^2/2 + 2 \pair{\nabla \ell, \nabla v} \sigma v.
    \end{align}
The last term on the right-hand side of (\ref{A_squared})
can be rewritten as
    \begin{align}\label{A_squared_cross}
2 \pair{\nabla \ell, \nabla v} \sigma v
= 
\pair{\sigma \nabla \ell, \nabla v^2}
= \div(v^2 \sigma \nabla \ell) - \div(\sigma \nabla \ell) v^2.
    \end{align}

Let us now study
    \begin{align}\label{SA}
SA = A \Box v + Aqv
= -\div (A \nabla v) + \pair{\nabla A, \nabla v} + Aqv.
    \end{align}
Here,
    \begin{align}\label{A_box_aux}
\pair{\nabla A, \nabla v}
= 2 \pair{\nabla \pair{\nabla \ell, \nabla v}, \nabla v} + \pair{\nabla (\sigma v), \nabla v}.
    \end{align}
We begin with the first term on the right-hand side of (\ref{A_box_aux}). Viewing $\nabla v$ as a vector field acting on a function,
    \begin{align*}
\pair{\nabla \pair{\nabla \ell, \nabla v}, \nabla v}
= \Hess \ell(\nabla v, \nabla v) + \Hess v(\nabla v, \nabla \ell).
    \end{align*}
On the other hand,
    \begin{align*}
2 \Hess v(\nabla v, \nabla \ell) 
= \div(\pair{\nabla v, \nabla v} \nabla \ell) + \Box \ell\;\pair{\nabla v, \nabla v}.
    \end{align*}
We turn now to the second term on the right-hand side of (\ref{A_box_aux}),
    \begin{align*}
\pair{\nabla (\sigma v), \nabla v}
&= 
\sigma \pair{\nabla v, \nabla v}
+ \frac12 \pair{\nabla \sigma, \nabla v^2}
\\&= 
\sigma \pair{\nabla v, \nabla v}
+ \frac12 \div(v^2\nabla \sigma)
+ \frac12 (\Box \sigma) v^2.
    \end{align*}
Combining the above gives
    \begin{align}\label{A_box}
A \Box v 
&= 
\tilde \sigma \pair{\nabla v, \nabla v} + 2 \Hess \ell(\nabla v, \nabla v) \\\notag&\qquad
+ \frac12 (\Box \sigma) v^2 + \div \left(
-A\nabla v + \pair{\nabla v, \nabla v} \nabla \ell + \frac12 v^2\nabla \sigma\right).
    \end{align}

Let us now consider the term
    \begin{align*}
Aqv = 2\pair{\nabla \ell, \nabla v} q v + \sigma q v^2,
    \end{align*}
where the first term on the right-hand side satisfies,
    \begin{align*}
2\pair{\nabla \ell, \nabla v} q v
= 
\pair{q \nabla \ell, \nabla v^2}
= \div(v^2 q \nabla \ell) - \div(q \nabla \ell)\, v^2.
    \end{align*}
Recall that $q= -\pair{\nabla \ell, \nabla \ell} - \tilde \sigma$ and that
    \begin{align*}
\div(\pair{\nabla \ell, \nabla \ell} \nabla \ell) 
&= 
-\pair{\nabla \ell, \nabla \ell} \Box \ell 
+ \pair{\nabla \ell, \nabla \pair{\nabla \ell, \nabla \ell}}
\\&=
-\pair{\nabla \ell, \nabla \ell} \Box \ell 
+ 2 \Hess \ell(\nabla \ell, \nabla \ell).
    \end{align*}
Therefore, by combining the above identities we obtain
    \begin{align}\label{A_q}
Aqv &= \left(-\tilde \sigma \pair{\nabla \ell, \nabla \ell}
+ 2 \Hess \ell(\nabla \ell, \nabla \ell)\right)v^2
\\\notag&\qquad
+ \div(\tilde \sigma \nabla \ell)\, v^2 - \sigma \tilde \sigma v^2 + \div(q v^2  \nabla \ell).
    \end{align}
The claim follows by combining (\ref{conj_box_square})--(\ref{SA}), (\ref{A_box}) and (\ref{A_q}).
\end{proof}

\subsection{Proof of Theorem~\ref{thm_uc}}

We assume that the Lorentzian manifold $(M,g)$, the point $p$ and the distribution $u$ are as in Theorem~\ref{thm_uc}. We begin by fixing a small number $\rho>0$ and define the set 
\begin{equation}
\Omega_\rho = \{q \in \mathscr E_p\,:\, r_p(q)>\rho\},
\end{equation}
where $r_p$ is the distance function defined by \eqref{def_r}. Recall from Lemma~\ref{lem_smoothness} that $u \in C^{\infty}(\xbar{\Omega_\rho})$. Recall also that $r_p$ is smooth on $\Omega_\rho$. 

For each $\epsilon\in (0,\rho)$, we let $F_{\rho,\epsilon}:(\rho,\infty)\to \R $ be defined by
$$ F_{\rho,\epsilon}(t)=4\log(t-\rho)+\tau_{\rho,\epsilon}\, (t-\rho)^2,$$
where $\tau_{\rho,\epsilon}>\epsilon^{-1}$ is a large constant depending on $\epsilon$ and $\rho$ that will be fixed later in the proof (see \eqref{tau}). Next, we consider the smooth weight function $\ell_{\rho,\epsilon}$ on $\Omega_\rho$ that is defined by
\begin{equation}\label{weight_fcn}\ell_{\rho,\epsilon}(q)= (F_{\rho,\epsilon}\circ r_p)(q)= 4\log (r_p(q)-\rho) + \tau_{\rho,\epsilon}\,(r_{p}(q)-\rho)^2, \quad \forall\,q\in \Omega_{\rho}.\end{equation}

We will apply Lemma~\ref{lem_basic_carleman} twice on the set 
\begin{equation}\label{U_epsilon}U_{\rho,\epsilon}=\Omega_{\rho+2\sqrt{\epsilon}},\end{equation} 
with the weight function $\ell_{\rho,\epsilon}$ above and with the function $v$ equal to 
$$v_1=e^{\ell_{\rho,\epsilon}}\,\Re{u},\quad \text{or}\quad v_2=e^{\ell_{\rho,\epsilon}}\,\Im{u},$$
where $\Re u$, $\Im u$ denote the real and imaginary parts of $u$ respectively. By taking a limit as $\epsilon$ approaches zero, we will show that the function $u$ vanishes identically on the set $\Omega_\rho$. Finally, by letting $\rho$ approach zero, we conclude the proof of the theorem.

In the remainder of this section and for the sake of brevity, we use the abbreviated notations $r$, $F$, $\ell$, $\tau$, $U$ and $\psi$ in place of $r_p$, $F_{\rho,\epsilon}$, $\ell_{\rho,\epsilon}$, $\tau_{\rho,\epsilon}$, $U_{\rho,\epsilon}$ and $\psi_{K,p}$ respectively. We will also hide the evaluation at a point $q$, writing for example $F(r)$ in place of $F_{\rho,\epsilon}(r(q))$. 

We record that
    \begin{equation}\label{F_prop}
    \begin{aligned}
F'(r) &= 4(r-\rho)^{-1} + 2\tau(r-\rho),\\
F''(r) &= -4 (r-\rho)^{-2}+2\tau.
    \end{aligned}
    \end{equation}
Since $r-\rho >2\sqrt{\epsilon}$ on $U$ and $\tau>\epsilon^{-1}$, it follows that
\begin{align}
\label{F''} 
F''(r)> \tau\quad \text{on $U$}.
\end{align}
We will also record the following simple calculation: 
    \begin{align}\label{Hess_ell}
\Hess \ell (X, X)
= 
F''(r) \pair{\nabla r, X}^2
+ F'(r) \Hess r(X, X).
    \end{align}
Using Lemma~\ref{lem_r} it follows that
    \begin{align}\label{Hess_ell2}
\Hess \ell (\nabla \ell, \nabla \ell)&=F'(r)^2 (F''(r) \pair{\nabla r, \nabla r}^2
+ F'(r) \Hess r(\nabla r, \nabla r))\\
\notag&=F'(r)^2\,F''(r).
    \end{align}

We fix a small $\rho>0$ and let $\epsilon \in (0,\rho)$. Recall from Lemma~\ref{lem_smoothness} that the function $v$ is smooth over $\xbar{\Omega_\rho}$ that contains $\xbar{U_{\rho,\epsilon}}$. Next, in view of the right hand side of Lemma~\ref{lem_basic_carleman}, we proceed to find a lower bound for the expression
$$ P+ 2\langle \nabla \ell ,\nabla v\rangle^2\quad \text{on $U$}.$$ 
Applying (\ref{Hess_ell}), we write
    \begin{multline*}
2\pair{\nabla \ell , \nabla v}^2 
+ \tilde \sigma \pair{\nabla v, \nabla v} + 2 \Hess \ell (\nabla v, \nabla v)\\
= 2F'(r)\Hess r(\nabla v,\nabla v) + \tilde \sigma \pair{\nabla v, \nabla v}
+ 2 ((F'(r)^2 + F''(r))\pair{\nabla r, \nabla v}^2\\
\geq (2\frac{F'(r)}{r}\psi+\tilde \sigma) \pair{\nabla v, \nabla v}+2(F'(r)^2 + F''(r)-\frac{F'(r)}{r}\psi)\pair{\nabla r, \nabla v}^2,
    \end{multline*}
where we are using Lemma~\ref{spacetime_conv_lem} in the last step, together with the fact that $F'(r(q))>0$ on $\Omega_\rho$. Next, we take 
$$\tilde \sigma(q)=-2\frac{F'(r(q))}{r(q)}\psi(q),\qquad \forall\, q \in \Omega_\rho.$$ 
Then, using \eqref{F_prop}--\eqref{F''} the previous bound reduces as follows: 
    \begin{align*}
&2\pair{\nabla \ell , \nabla v}^2 
+ \tilde \sigma \pair{\nabla v, \nabla v} + 2 \Hess \ell (\nabla v, \nabla v)
\\&\qquad \geq 2(F'(r)^2 + F''(r)-\frac{F'(r)}{r}\psi)\pair{\nabla r, \nabla v}^2
\\&\qquad \geq  2F'(r)(F'(r)-\frac{1}{r}\psi)\pair{\nabla r, \nabla v}^2
\\&\qquad \geq  2F'(r)(\epsilon^{-\frac{1}{2}}-\frac{1}{r}\psi)\pair{\nabla r, \nabla v}^2\geq 0, \quad \text{on $U$,}
\end{align*}
where the latter inequality holds for all $\epsilon$ sufficiently smaller than $\rho$, since $\frac{\psi}{r}$ is a bounded smooth function on $\xbar{\Omega_\rho}$. Using this bound together with \eqref{Hess_ell2} we obtain 
 \[
    \begin{aligned}
& P+2\langle\nabla \ell ,\nabla v\rangle^2\geq \left(-\tilde \sigma \pair{\nabla \ell , \nabla \ell } + 2 \Hess \ell (\nabla \ell , \nabla \ell )\right)v^2
\\&\qquad\geq
2 F'(r)^2\left( F'(r)\psi\,r^{-1}+F''(r) \right)v^2\geq 2F'(r)^2 F''(r) v^2 \quad \text{on $U$},
    \end{aligned}
\]
where in the last step we used the fact that $\psi$, $r$ and $F'(r)$ are positive functions on $U$. Thus,
 \begin{equation}\label{bound_key_1}
P+2\langle\nabla \ell ,\nabla v\rangle^2\geq 2\tau (r-\rho)^2 \lambda^2\,v^2,\quad \text{on $U$,}
\end{equation}
where 
\begin{equation}
\label{lambda_def}
\lambda=\lambda(r)=\tau+(r-\rho)^{-2}.
\end{equation}

Next, we consider the terms of the expression $R$ in Lemma~\ref{lem_basic_carleman}. We claim that the following bound holds:
\begin{equation}\label{bound_key_2} |R| \leq c_\rho\, \lambda^2v^2\quad \text{on $U$},\end{equation}
where $c_\rho>1$ is a constant only depending on $\rho$.

In order to prove \eqref{bound_key_2}, we will first record some bounds. First, note that
\begin{equation}
\label{F'_bound}
|F'(r)|+|F''(r)| \leq C\lambda,\quad \text{on $\Omega_\rho$}
\end{equation}
for some constant $C>0$ independent of $\rho$ and $\epsilon$. Next, note that
   \begin{equation}\label{Box_l}
\Box \ell = -\div(F'(r) \nabla r) 
= -F''(r) + F'(r) \Box r.
    \end{equation}
Recalling also that $\sigma=\tilde{\sigma}-\Box \ell$, it follows that
 \begin{equation}
 \label{bound_box_sigma}
 |\Box \ell|+|\sigma| \leq C_\rho \lambda,\quad \text{on $\Omega_\rho$},
 \end{equation}
 for some $C_\rho>0$ that depends only on $\rho$.
 
 We return to the claim \eqref{bound_key_2} and start with the term $\div((\Box \ell)\nabla\ell)v^2$ in $R$. There holds
    \begin{align*}
\div((\Box \ell)\nabla\ell)= \pair{\nabla \Box \ell, \nabla \ell} - (\Box \ell)^2.
    \end{align*}
To simplify notation we will write $F=F(r)$. Using \eqref{Box_l} we write
    \begin{align*}
\div((\Box \ell)\nabla\ell)=-F'''F'+F''F'\Box r+F'^2\langle\nabla \Box r,\nabla r \rangle-(\Box \ell)^2.
    \end{align*}
Next, using \eqref{F_prop} and noting that $r$ is a smooth bounded function on $U$, it is clear that the claimed estimate \eqref{bound_key_2} holds for the first term in $R$ with $c_\rho$ sufficiently large depending on $\rho$. We move on to analyze the term $(-\sigma\tilde{\sigma}+\frac{1}{2}\sigma^2)v^2$ in the expression for $R$. There holds
    \begin{align*}
-\sigma\tilde{\sigma}+\frac{1}{2}\sigma^2=-\frac{1}{2}\tilde{\sigma}^2+\frac{1}{2}(\Box \ell)^2= -\frac{1}{2} F'^2 \frac{\psi^2}{r^2}+\frac{1}{2}(-F''+F'\Box r)^2.
    \end{align*}
Again, it is clear that the claimed bound \eqref{bound_key_2} holds for this term with $c_\rho$ sufficiently large depending on $\rho$. It remains to analyze the term $(\Box \sigma)v^2$:
    \begin{align*}
\Box \sigma &= \Box \tilde \sigma - \Box( \Box \ell)\\
&= -F''''-2F' \Box(\frac{\psi}{r})+2(F''' -F''\Box r)\frac{\psi}{r}-F' \Box(\Box r)+2F'''\Box r\\
&\quad\,+4 F'' \langle \nabla r,\nabla (\frac{\psi}{r})\rangle-F''(\Box r)^2+2F''\langle \nabla r,\nabla\Box r\rangle.
\end{align*}
The claim follows again by using the expression \eqref{F_prop} and choosing $c_\rho$ sufficiently large depending on $\rho$. This completes the proof of the bound \eqref{bound_key_2}.

We are ready to fix the choice of $\tau>\epsilon^{-1}$. We set,
\begin{equation}\label{tau}\tau= a_\rho\, \epsilon^{-1},\end{equation}
where $a_\rho= \max\{\|V\|^2_{L^{\infty}(M)},c_\rho\}$ and $c_\rho$ is as in \eqref{bound_key_2}. Using the bounds \eqref{bound_key_1} and \eqref{bound_key_2} together with the fact that $r-\rho> 2\sqrt{\epsilon}$ on $U$, it follows that
\begin{equation}
\label{combined_pos}
P+2\langle\nabla \ell,\nabla v\rangle^2+R \geq  2a_\rho v^2,\quad \text{on $U$}.
\end{equation}

Next, writing $u_1= \Re u$, $u_2= \Im u$, and recalling that $(\Box+V)(u_1+iu_2)=0$ and that $v_j=e^{\ell}u_j$ for $j=1,2$, we obtain
$$ \frac{1}{2}|V (v_1+iv_2)|^2 = \frac{1}{2}|e^{\ell}\Box(e^{-\ell} v_1)|^2+ \frac{1}{2}|e^{\ell}\Box(e^{-\ell} v_2)|^2.$$
Integrating the latter expression over $U$, applying Lemma~\ref{lem_basic_carleman} with $v=v_1$ and $v=v_2$ and finally using the estimate \eqref{combined_pos}, we deduce that
$$2a_\rho \int_{U}e^{2\ell} |u|^2\,dV_g\leq \left|\int_{U} \div B_1 \,dV_g \right|+\left|\int_{U} \div B_2 \,dV_g \right|+\frac{1}{2}\int_U |V|^2e^{2\ell} |u|^2\,dV_g,$$
where $B_j$, $j=1,2$ is as in Lemma~\ref{lem_basic_carleman} with $v=v_j$.
Since $a_\rho> \|V\|^2_{L^{\infty}(M)}$, this reduces to
\begin{equation}\label{carleman_eq}
a_\rho \int_{U}e^{2\ell} |u|^2\,dV_g\leq \left|\int_{U} \div B_1 \,dV_g \right|+\left|\int_{U} \div B_2 \,dV_g \right|.
\end{equation}
Using the divergence theorem together with the fact that both $u$ and $\p_\nu u$ vanish identically on the set $\mathscr E_p\cap \Sigma$, we write for $j=1,2,$
\begin{equation}\label{div_thm} \int_U \div B_j \,dV_g = \int_\Gamma \langle B_j,\nu\rangle\,dV_g,\end{equation}
where $\Gamma=\Gamma_{\rho,\epsilon}= \{ r-\rho=2\sqrt{\epsilon}\} \cap \xbar{\Omega_\rho}$ and $\nu$ is the unit spacelike normal vector to $\Gamma$. Next, noting that 
\begin{equation}
\label{r_boundary}
(r-\rho)|_{\Gamma}=2\epsilon^{\frac{1}{2}},\end{equation} 
we write
\begin{equation}\label{exp_on_boundary} e^{2\ell} |_{\Gamma} = (r-\rho)^8 e^{2\tau(r-\rho)^2}|_{\Gamma} = 2^8\epsilon^4 e^{8a_\rho}.\end{equation}
We claim that the following bound holds: 
\begin{equation}\label{B_bound} \left|\langle B_j,\nu\rangle\right| \leq C_\rho \lambda^3 e^{2\ell}(u_j^2 + |\nabla u_j|^2), \quad \text{pointwise on $\Gamma$},\end{equation}
where $|\nabla u_j|$ is the norm of $\nabla u_j$ with respect to some auxiliary Riemannian metric and $C_\rho$ is a generic positive constant only depending on $\rho$. 

To this end, let us first observe that
    \begin{align*}
|\nabla v_j| 
= |e^\ell (u_j \nabla \ell + \nabla u_j)|
= |e^\ell (F' u_j \nabla r + \nabla u_j)|
\leq C_{\rho} e^\ell F'\,(|u_j|+|\nabla u_j|),
    \end{align*}
pointwise on $\Gamma$, for a positive constant $C_{\rho}$ only depending on $\rho$. Using this expression and recalling the bound \eqref{bound_box_sigma} we obtain the following bounds for the terms in the definition of $B$:
    \begin{align*}
|\pair{\nabla \ell, \nabla v_j}\nabla v_j|+|\sigma v_j \nabla v_j|+|\pair{\nabla v_j, \nabla v_j} \nabla \ell|
\leq C_{\rho} e^{2\ell} \lambda^3(|u_j|^2+|\nabla u_j|^2), 
    \end{align*} 
    \begin{align*}
|\Box \ell v_j^2 \nabla \ell|+|\pair{\nabla \ell, \nabla \ell} v_j^2 \nabla \ell|+|v_j^2 \nabla \sigma|
\leq C_{\rho} e^{2\ell} \lambda^3|u_j|^2,
    \end{align*}
for some constant $C_{\rho}$ that only depends on $\rho$. This completes proof of the bound \eqref{B_bound}. Combining \eqref{B_bound} with equations \eqref{tau}, \eqref{r_boundary} and \eqref{exp_on_boundary} we write
$$  \left| \int_{\Gamma} \langle B_j,\nu\rangle\,dV_g \right|\leq C_{\rho} \int_{\Gamma}((r-\rho)^{-2}+\tau)^3 e^{2\ell}(u_j^2 + |\nabla u_j|^2)\,dV_g \leq C_\rho\, \epsilon,$$
for some positive constant $C_\rho$ only depending on $\rho$. 

Next, combining the latter estimate with the Carleman estimate \eqref{carleman_eq} and equation \eqref{div_thm} it follows that
$$ \int_{U} (r-\rho)^{8}\,|u|^2\,dV_g \leq C_\rho \epsilon.$$
Therefore, taking the limit as $\epsilon$ tends to zero and noting that $U=U_{\rho,\epsilon} \to \Omega_\rho$ as $\epsilon$ approaches zero, we conclude that
$$ \int_{\Omega_\rho}(r-\rho)^{8}\,|u|^2\,dV_g=0.$$
Hence $u=0$ on the set $\Omega_\rho$. Finally, letting $\rho \to 0$, we conclude that $u=0$ on $\mathscr E_p$ as claimed.

\section{Exact controllability in rough Sobolev spaces}
\label{ctrl_sec}
Our main goal in this section is to prove the following proposition.
\begin{proposition}[Exact controllability]
\label{ctrl_prop}
Let $(M,g)$, $T_0$ and $p_0$ be as in Theorem~\ref{thm_main}. Let $s \geq 0$, $V \in C^{\infty}(M)$ and suppose that $w_0$ and $w_1$ are compactly supported distributions in the interior of $M_0$ such that
$$(w_0,w_1) \in H^{-s}(M_0)\times H^{-s-1}(M_0).$$  Then, given any $T_1 \in [T_0,T]$, there exists a compactly supported distribution
$$f \in H^{-s}(\Sigma)\cap \mathcal E'((-T,T_0)\times \p M_0),$$
such that the solution $u\in H^{-s}(M)$ to the equation \eqref{pf0} with boundary value $f$ satisfies
 $$(u|_{t=T_1},\p_t u|_{t=T_1})= (w_0,w_1)\quad \text{on $M_0$.}$$
\end{proposition}

\subsection{Direct problem in $H^{-s}$-spaces with $s\geq 0$}

We first study the direct problem \eqref{pf0} with data in smooth Sobolev spaces, and then proceed to study rougher data via transposition. 
\begin{lemma}
\label{lem_smooth}
Let $s\geq 0$ and $V \in C^{\infty}(M)$. Consider the initial boundary value problem
\begin{equation}\label{pf_full}
\begin{aligned}
\begin{cases}
\Box u+Vu=F\,\quad &\text{on $M$},
\\
u=f\,\quad &\text{on $\Sigma$,}\\
(u|_{t=T},\p_{t}u|_{t=T})=(\phi_0,\phi_1) \,\quad &\text{on $M_0$.}
\end{cases}
 \end{aligned}
\end{equation}
Given any $F \in H^{s}_0(M)$, $u_0 \in H^{s+1}_0(M_0)$, $u_1 \in H^s_0(M_0)$ and $f\in H^{s+1}_0(\Sigma)$, there exists a unique solution $u$ in the energy space
\begin{equation}\label{smooth_sol}
H^{s+1}(M)\cap C^1(-T,T;H^{s}(M_0))\cap  C(-T,T;H^{s+1}(M_0)),
\end{equation}
and the dependence on the data is continuous. Moreover, $\p_\nu u\in H^{s}(\Sigma)$. 
\end{lemma} 
This lemma follows from \cite[Theorem 4.1]{LLT} by using standard techniques of obtaining higher regularity for solutions to the wave equation, see for example \cite[Remark 2.10]{LLT} or \cite[Chapter 7, Theorem 6]{Evans}. Next, we state a lemma on solving the wave equation in rough Sobolev scales via a standard transposition argument. For the convenience of the reader we have included the proof. 

\begin{lemma}
\label{rough_direct}
Let $s\geq 0$, $V\in C^{\infty}(M)$ and let 
$$f  \in H^{-s}(\Sigma)\cap \mathcal E'((-T,T_0)\times \p M_0),$$
for some $|T_0|< T$. Then, there exists a unique solution $u \in H^{-s}(M)$ to the equation
\begin{equation}\label{pf_rough}
\begin{aligned}
\begin{cases}
\Box u+Vu=0\,\quad &\text{on $M$},
\\
u=f\,\quad &\text{on $\Sigma$,}\\
(u|_{t=-T},\p_{t}u|_{t=-T})=0 \,\quad &\text{on $M_0$.}
\end{cases}
    \end{aligned}
\end{equation}
There holds
$$\p_\nu u|_{\Sigma}\in H^{-s-1}(\Sigma).$$ 
Moreover, for any $t_0\in [T_0,T]$, there holds
$$(u|_{t=t_0},\p_t u|_{t=t_0}) \in H^{-s}(M_0)\times H^{-s-1}(M_0),$$
\end{lemma}

\begin{proof}
Uniqueness follows by using propagation of singularities \cite[Theorem 24.5.3]{Ho3} and the fact $u=0$ is the only smooth solution to \eqref{pf_rough} subject to $f=0$. 

To show existence of a solution, we suppose for the moment that $f \in C^{\infty}_c((-T,T_0)\times M_0)$ and let $u$ be the solution to \eqref{pf_rough} with the boundary value $f$. Next, let $u^F$ denote the unique solution to \eqref{pf_full} with source $F$ and zero boundary and final data. We have
$$ \int_M u(t,x) \, F(t,x)\,dV_g = -\int_\Sigma \p_\nu u^F(t,x) \,f(t,x)\,dV_g.$$
Thus, the transpose of the map $F \mapsto \p_\nu u^F|_{\Sigma}$ is $f \mapsto u$. By Lemma~\ref{pf_full}, the former map is continuous from $H^{s}_0(M)$ to $H^s(\Sigma)$. Hence, the latter is continuous from $H^{-s}(\Sigma)\cap \mathcal E'(\widetilde\Sigma)$ to $H^{-s}(M)$ for any compact set $\widetilde\Sigma$ in $(-T,T_0)\times \p M_0$. 

To show the claim about the trace on $\Sigma$, we let $f$ and $u$ be as above and given any $h \in H^{s+1}_0(\Sigma)$, we write $u^h$ for the unique solution to \eqref{pf_full} subject to the boundary value $h$ and zero source and final data. Then,
$$ \int_\Sigma \p_\nu u^h(t,x)\,f(t,x)\,dV_g=\int_\Sigma \p_\nu u(t,x)h(t,x)\,dV_g.$$
Hence, the transpose of the map $h \mapsto \p_\nu u^h|_{\Sigma}$ is $f \mapsto \p_\nu u|_{\Sigma}$. By Lemma~\ref{pf_full}, the former map is continuous from $H^{s+1}_0(\Sigma)$ to $H^s(\Sigma)$. Hence, the latter is continuous from $H^{-s}(\Sigma)\cap \mathcal E'(\widetilde\Sigma)$ to $H^{-s-1}(\Sigma)$ for any compact set $\widetilde\Sigma$ in $(-T,T_0)\times \p M_0$. 

Finally, we consider the last claim and only consider the case $t_0=T$. The case $t_0\in [T_0,T)$ is analogous. Let $\phi=(\phi_0,\phi_1) \in H^{s+1}_0(M_0)\times H^{s}(M_0)$ and denote by 
$u^\phi\in H^{s+1}(M),$ the unique solution to \eqref{pf_full} with $F=f=0$ and final data $\phi$. We have
$$ \int_{\Sigma}f(t,x) \p_\nu u^\phi(t,x)\,dV_g = \int_{M_0} u(T,x)\phi_1(x)\,dV_{g} -\int_{M_0} \p_t u(T,x)\phi_0(x)\,dV_{g}.$$
Thus, the transpose of the map $\phi \mapsto \p_\nu u^\phi|_{\Sigma}$ is $f \mapsto (u|_{t=T},\p_t u|_{t=T}).$ By Lemma~\ref{pf_full} the former map is continuous from $H^{s+1}_0(M_0)\times H^{s}(M_0)$ to $H^s(\Sigma)$. Hence, the latter is continuous from $H^{-s}(\Sigma)\cap \mathcal E'(\widetilde\Sigma)$ to $H^{-s}(M_0)\times H^{-s-1}(M_0)$ for any compact set $\widetilde\Sigma$ in $(-T,T_0)\times \p M_0$. 
\end{proof}

\subsection{Proof of exact controllability}

The aim of this section is to prove Proposition~\ref{ctrl_prop}. The proof is based on the following observability estimate.
\begin{lemma}[Observability estimate]
\label{lem_obs}
Let $(M,g)$, $T_0$ and $p_0$ be as in Theorem~\ref{thm_main}. Let $V \in C^{\infty}(M)$, let $\delta>0$ be sufficiently small and let $\chi\in C^{\infty}_c((-T,T_0)\times \p M_0)$ satisfy $\chi=1$ on $(-T+\delta,T_0-\delta)\times \p M_0$. Then, given any $s \geq 0$ there holds:
\begin{equation}
\label{obs}
\|\phi_0\|_{H^{s+1}_0(M_0)}+\|\phi_1\|_{H^s_0(M_0)} \leq C_{\textrm{obs}}\,\|\chi u\|_{H^s_0(\Sigma)},
\end{equation}
where $u$ denotes the unique solution to the wave equation \eqref{pf_full} on $M$, subject to the source $F=0$, boundary value $f=0$ and final data $$(\phi_0,\phi_1)\in H^{s+1}_0(M_0)\times H^s_0(M_0).$$
\end{lemma}

Before presenting the proof we note that in the following discussion the term {\em compressed generalized bicharacteristic} is as defined in \cite{BLR}. We write $\iota : \Sigma \to M$ for the natural inclusion and write $\iota^* : T^* M \to T^* \Sigma$ by its pullback. We also recall from \cite{BLR} that a pair $(y,\eta)\in T^*\Sigma$ is said to be {\em nondiffractive} if the preimage $(\iota^*)^{-1}(y,\eta)$ contains two distinct null vectors, or if it contains a unique null vector $(y, \tilde \eta)$ and for all $s_0>0$ there is $|s|<s_0$ such that $\gamma_{y,\tilde \eta}(s) \notin M$, when $M$ is extended to a slightly larger manifold without boundary.

\begin{proof}
We begin by choosing $\delta>0$ to be sufficiently small so that $$\xbar{\mathscr E_{p_0}}\subset (-T+\delta,T_0-\delta)\times M_0.$$ 
We claim that that every bicharacteristic $(\gamma_{q,\xi}(s),\dot{\gamma}_{q,\xi}(s))$ with $(q,\xi)\in T^*\mathscr E_{p_0}$ must hit $T^*\Gamma$ in a nondiffractive point, where 
$$\Gamma=(-T+\delta,T_0-\delta)\times \p M_0.$$ 
Here and in the remainder of this proof we are identifying vectors and covectors using $g$ and use the same notation for both of them. We also write $[a,b]$ for the maximal interval of $\gamma_{q,\xi}$ in $M$.

By Lemma~\ref{lem_cutlocus}, $\gamma_{v}$ can not intersect $J^+(p_0)\cup J^-(p_0)$ twice. Thus $\gamma_{q,\xi}(s)\in \mathcal E_{p_0}$ for all $s \in [0,b]$ or for all $s \in [a,0]$. We consider only the former case, the other case being analogous.
We write 
    \begin{align*}
(y,\eta) =\iota^* (\gamma_{q,\xi}(b), \dot \gamma_{q,\xi}(b)) \in T^* \Sigma.
    \end{align*}
It is enough to show that $(y,\eta)$ is nondiffractive. To make the discussion explicit, let us consider boundary normal coordinates $(x',r) \in \Sigma \times (-r_0,r_0)$ such that 
    \begin{align*}
g(x',r) = \begin{pmatrix}
h(x',r) & 0 \\
0 & 1
\end{pmatrix}
    \end{align*}
and $r \geq 0$ for $(x',r) \in M$. In these coordinates, $\dot \gamma_{q,\xi}(b) = (\eta, \rho)$ for some $\rho \leq 0$.

Case $\rho < 0$. As $(\eta, \rho)$ is a null vector, so is $(\eta, -\rho)$. Thus $(\iota^*)^{-1}(y,\eta)$ contains two distinct null vectors.

Case $\rho = 0$. Then $(y, \tilde \eta) = (y, 0; \eta, 0)$ is the only point in 
$(\iota^*)^{-1}(y,\eta)$.
But $(\eta, 0) = \dot \gamma_{q,\xi}(b)$ and 
there are $s > b$ arbitrarily close to $b$ such that
$\gamma_{q,\xi}(s) \notin M$, due to maximality of $[a,b]$. 
 
Combining with the fact that $M$ is connected, it follows that every compressed generalized bicharacteristic in $T^*M$ must hit $T^*\Gamma$ in a nondiffractive point. Thus, the claim follows by invoking \cite[Theorem 3.3--Corollary 3.7]{BLR}, if given any $v \in C^{\infty}(M)$ satisfying the equation $(\Box+V)v=0$ on $M$ and $v|_{\Sigma}=\p_\nu v|_{\Gamma}=0$, there holds $v=0$ everywhere on $M$.

To verify this claim, we write $p_0=(t_0,x_0)$ and note that since $\xbar{\mathscr E_{p_0}}\subset (-T+\delta,T_0-\delta)\times M_0$ and since $v|_{\Gamma}=\p_\nu v|_{\Gamma}=0$, we can apply Theorem~\ref{thm_uc} to deduce that $(v|_{t=t_0},\p_t v|_{t=t_0})=0$ on $M_0\setminus \{x_0\}$. By smoothness of $v$ in $M$ it follows that $(v|_{t=t_0},\p_t v|_{t=t_0})=0$ on $M_0$. Since $(\Box+V)v=0$ on $M$ and $v|_{\Sigma}=0$ we conclude that $v=0$ everywhere on $M$.
\end{proof}

It is well-known that observability estimates imply exact controllability. We will give a proof for the convenience of the reader.

\begin{proof}[Proof of Proposition~\ref{ctrl_prop}]
We choose $\delta>0$ sufficiently small as in the proof of Lemma~\ref{lem_obs} and define the smooth non-negative function $$\chi\in C^{\infty}_c((-T,T_0)\times \p M_0)$$ such that $\chi=1$ on $(-T+\delta,T_0-\delta)\times \p M_0$. Given any $F \in H^{-s}(\Sigma)$, we write $u \in H^{-s}(M)$ for the solution to \eqref{pf_rough} with the boundary value $f=\chi F$. Thus, we can view $\chi F$ as an element in 
\begin{equation}\label{right_space}
H^{-s}(\Sigma)\cap \mathcal E'((-T,T_0)\times \p M_0).\end{equation} 
By Lemma~\ref{rough_direct} we obtain that 
$$ (u|_{t=T_1},\p_t u|_{t=T_1})\in H^{-s}(M_0)\times H^{-s-1}(M_0).$$
Given any $\phi=(\phi_0,\phi_1) \in H^{s+1}_0(M_0)\times H^{s}_0(M_0)$, let $u^\phi \in H^{s+1}(M)$ be the solution to
 \begin{equation}\label{pf_reverse}
\begin{aligned}
\begin{cases}
\Box u^\phi+Vu^\phi=0\,\quad &\text{on $M$},
\\
u^\phi=0\,\quad &\text{on $\Sigma$,}\\
(u^\phi(T_1,x),\p_{t}u^\phi(T_1,x))=(\phi_0,\phi_1) \,\quad &\text{on $M_0$.}
\end{cases}
    \end{aligned}
\end{equation}
We define the continuous linear map $$\mathcal T: H^{s+1}_0(M_0)\times H^{s}_0(M_0)\to H^{s}_0(\Sigma),$$
through $\mathcal T \phi =\chi \p_\nu u^\phi\big|_{\Sigma}$. Note that $\mathcal T\phi$ is compactly supported on the set $(-T,T_0)\times \p M_0$. We write
$$ (\mathcal T \phi, F)_{(-T,T_1)\times \p M_0} = (\phi_0,\p_t u)_{\{T_1\}\times M_0}-(\phi_1,u)_{\{T_1\}\times M_0},$$
where the pairings $(\cdot,\cdot)_{(-T,T_1)\times \p M_0}$ and $(\cdot,\cdot)_{\{T_1\}\times M_0}$ are as defined in Section~\ref{prelim_sec_not}, that is, the pairing on the left hand side is in the generalized sense $$H^s_0((-T,T_1)\times \p M_0)\times H^{-s}((-T,T_1)\times \p M_0),$$ while the ones on the right hand side are in the generalized sense $$H^{s+1}_0(M_0)\times H^{-s-1}(M_0)\quad \text{and}\quad H^{s}_0(M_0)\times H^{-s}(M_0)$$ respectively. It follows from the latter identity that the transpose
$$\mathcal T^*: H^{-s}(\Sigma)\to H^{-s}(M_0)\times H^{-s-1}(M_0),$$   
defined by $$\mathcal T^*F = (u|_{\{T_1\}\times M_0},\p_t u|_{\{T_1\}\times M_0}),$$ is a continuous linear map. To conclude the proof, it suffices to show that $\mathcal T^*$ is surjective. This is true, since $\mathcal T$ is injective and has closed range, thanks to Lemma~\ref{lem_obs}. We emphasize that since $F \in H^{-s}(\Sigma)$ and since $\chi$ is compactly supported on $(-T,T_0)\times \p M_0$, it follows that the boundary data $f= \chi F$ that drives the solution $u$ to the target state $(w_0,w_1)$ is indeed in the right space \eqref{right_space}. 
\end{proof}

\section{Proof of Theorem~\ref{thm_main}}
\label{thm_sec}

This section is concerned with the proof of Theorem~\ref{thm_main}. We start with a proposition.
\begin{proposition}
\label{prop_final}
Let the hypotheses of Theorem~\ref{thm_main} be satisfied. Let $p=(T_1,x_0)$ for some $T_1 \in (T_0,T)$ and $x_0\in M_0^{\textrm{int}}$ be such that $\xbar{\mathscr E_p}\subset (T_0,T)\times M_0$. Then, if 
$$ \Lambda_{V_1}= \Lambda_{V_2},\quad \text{for some $V_1,V_2 \in C^{\infty}(M)$},$$
it follows that given any $f \in H^{\frac{n+1}{2}}_0((-T,T_1)\times \p M_0)$, there holds:
$$u_f^{(1)}(p)=\omega\,u_f^{(2)}(p)\quad \text{for some constant $\omega$ that is independent of $f$},$$ 
where for each $j=1,2$, the notation $u_f^{(j)}$ stands for the unique solution to \eqref{pf0} subject to the potential $V_{j}$ and the Dirichlet data $f$. 
\end{proposition}

Let us emphasize that we are identifying $H^{\frac{n+1}{2}}_0((-T,T_1)\times \p M_0)$ as a subspace of $H^{\frac{n+1}{2}}_0(\Sigma)$ by extending its elements by zero on the set $(T_1,T)\times \p M_0$. In the following proof, in order to avoid confusion, we will systematically use the notation $h$ for rougher Dirichlet data on the boundary, while $f$ is reserved for smoother Dirichlet data on the boundary as in the statement of the proposition.

\begin{proof}[Proof of Proposition~\ref{prop_final}]
We recall that by Lemma~\ref{pf_rough}, the equality $\Lambda_{V_1}=\Lambda_{V_2}$ on $H^1_0(\Sigma)$ extends by continuity to $H^{-s}(\Sigma)$ for all $s\geq 0$. Let $p=(T_1,x_0)$ be as in the statement of the proposition and consider $$\delta_{x_0}\in H^{-\frac{n+1}{2}}(M_0)$$ defined by
$$(\delta_{x_0}, v)_{\{T_1\}\times M_0}=v(T_1,x_0)\quad \text{for all $v \in C^{\infty}_c(\{T_1\}\times M_0^{\textrm{int}})$},$$
where the pairing above is with respect to the natural volume form on the hypersurface $\{T_1\}\times M_0$, as defined in Section~\ref{prelim_sec_not}. Applying Proposition~\ref{ctrl_prop}, we deduce that there exists a distribution $h$ in the space
\begin{equation}\label{rough_space} H^{-\frac{n-1}{2}}(\Sigma)\cap \mathcal E'((-T,T_0)\times \p M_0),\end{equation} 
such that $u_h^{(1)}$ satisfies 
\[(u_h^{(1)}|_{t=T_1},\p_t u_h^{(1)}|_{t=T_1})=(0,\delta_{x_0}), \quad \text{on $M_0$},\]
where we recall that the traces are well defined by Lemma~\ref{rough_direct}. 

Next, we denote by $\mathscr T^{(1)}$ the set of distributions $h$ in the space \eqref{rough_space} such that the solution $u_h^{(1)}$ satisfies 
\begin{equation}\label{def_T1} (u_h^{(1)}|_{t=T_1},\p_t u_h^{(1)}|_{t=T_1})\in \linspan{\{(0,\delta_{x_0})\}}, \quad \text{on $M_0$}.\end{equation}
We also define the set $\mathscr T^{(2)}$ analogously with $V_1$ replaced by $V_2$. We claim that the two sets $\mathscr T^{(1)}$ and $\mathscr T^{(2)}$ are identical. 

We show the inclusion $\mathscr T^{(1)} \subset \mathscr T^{(2)}$, and the opposite inclusion follows by symmetry. Let $h \in \mathscr T^{(1)}$ and recall that equation \eqref{def_T1} holds. Since $h$ belongs to the space \eqref{rough_space}, it follows from Lemma~\ref{rough_direct} that $u_h^{(1)}\in H^{-\frac{n-1}{2}}(M)$. Noting that $\xbar{\mathscr E_p}$ is disjoint from the set $\{T_0,T\}\times M_0$ and that $h$ is compactly supported in $(-T,T_0)\times \p M_0$, it follows that $u_h^{(1)}|_{\Sigma}$ is vanishing on the open set $\mathscr E_p \cap \Sigma$. Since $(u_h^{(1)}|_{t=T_1}, \p_t u_h^{(1)}|_{t=T_1})$ is supported at the point $x_0$ and since $u_h^{(1)}|_{\Sigma}$ is vanishing on $\Sigma \cap \mathscr E_p$, it follows by the finite speed of propagation for wave equation that $\p_\nu u_h^{(1)}|_{\Sigma}$ must also vanish on the set $\mathscr E_p\cap \Sigma$. 

Since $\Lambda_{V_1}h=\Lambda_{V_2}h$, it follows that the Cauchy data of the distribution $u_h^{(2)}\in H^{-\frac{n-1}{2}}(M)$ must vanish on the set $\mathscr E_p\cap \Sigma$. Applying our unique continuation result, Theorem~\ref{thm_uc}, it follows that $u_h^{(2)}|_{\mathscr E_p}=0$. Thus, we conclude that the traces $(u_h^{(2)}|_{t=T_1}, \p_t u_h^{(2)}|_{t=T_1}) \in H^{-\frac{n-1}{2}}(M_0)\times  H^{-\frac{n+1}{2}}(M_0)$ are supported at the point $\{x_0\}$. Consequently, the traces must be equal to a linear combination of $\delta_{x_0}$ and its derivatives on $M_0$. Since $\delta_{x_0} \in H^{-\frac{n+1}{2}}(M_0)$ while its derivatives are not in this space, it follows that $u_h^{(2)}|_{t=T_1}=0$ and that $\p_t u_h^{(2)}|_{t=T_1} \in \linspan{\{\delta_{x_0}\}}$. 

We have shown that $\mathscr T^{(1)}=\mathscr T^{(2)}$. Let us now consider $$f \in H^{\frac{n+1}{2}}_0((-T,T_1)\times \p M_0)$$ and let $h \in \mathscr T^{(1)}$ be such that 
$$(u_h^{(1)}|_{t=T_1},\p_t u_h^{(1)}|_{t=T_1})= (0,\delta_{x_0})\quad \text{on $M_0$}.$$
Recall that such $h$ exists thanks to Proposition~\ref{ctrl_prop}. In view of the equality $\mathscr T^{(1)}=\mathscr T^{(2)}$, we must also have:
\begin{equation}
\label{eq_h_1}
(u_h^{(2)}|_{t=T_1},\p_t u_h^{(2)}|_{t=T_1})=(0,\omega\,\delta_{x_0}),\quad \text{on $M_0$},
\end{equation}
for some constant $\omega$ that is independent of $f$.

Using Lemma~\ref{lem_smooth}, Lemma~\ref{rough_direct} and recalling that the distribution $h$ is compactly supported in $(-T,T_0)\times \p M_0$, we write, for each $j=1,2$,
    \begin{equation}
    \label{eq_h_2}
    \begin{aligned}
&(\Lambda_{V_j}h, f)_{(-T,T_1)\times\p M_0}-(h,\Lambda_{V_j}f)_{(-T,T_1)\times \p M_0}\\ 
&\qquad=(\p_t u_h^{(j)}, u_f^{(j)})_{\{T_1\}\times M_0}-(u_h^{(j)}, \p_t u_f^{(j)})_{\{T_1\}\times M_0}=c_{j} \,u_f^{(j)}(p),
    \end{aligned}
    \end{equation}
where $c_1=1$ and $c_2=\omega$ is the constant in \eqref{eq_h_1}. We emphasize here that the former two pairings are in the generalized sense $H^{-\frac{n+1}{2}}\times H^{\frac{n+1}{2}}_0$ and $H^{-\frac{n-1}{2}}\times H^{\frac{n-1}{2}}_0$ on $(-T,T_1)\times \p M_0$ while the latter two pairings are in the generalized sense $H^{-\frac{n+1}{2}}\times H^{\frac{n+1}{2}}_0$ and $H^{-\frac{n-1}{2}}\times H^{\frac{n-1}{2}}_0$ on $\{T_1\}\times M_0$.

Finally, we use equality of the Dirichlet-to-Neumann maps $\Lambda_{V_1}=\Lambda_{V_2}$ with \eqref{eq_h_2} to conclude that
$$ u_f^{(1)}(p)=\omega\, u_f^{(2)}(p).$$
\end{proof}
We are now ready to prove the main theorem, but first we state a lemma about point values of solutions to the wave equation that will be proved in Appendix~\ref{appendix A}. We recall that given $v=(p,\xi)\in L^+M$, with $L^+M$ denoting the bundle of future pointing null vectors on $M$, the notation $\gamma_v:I\to M$ stands for the inextendible null geodesic in $M$ with initial data $v$, as in \eqref{affine}.

\begin{lemma}\label{asymp_value}
Let $V\in C^{\infty}(M)$, $p=(T_1,x_0) \in M^{\textrm{int}}$, $v=(p,\xi)\in L^+M$ and suppose that the inextendible null geodesic $\gamma_v:[a,b]\to M$ satisfies $\gamma_v(a)\in \Sigma$. Let $\epsilon>0$ be small enough so that the null geodesic $\gamma_v$ is disjoint from $[T_1-\epsilon,T_1+\epsilon]\times \p M_0$. There exists a sequence $$\{f_{j}\}_{j=1}^{\infty}\subset C^{\infty}_c((-T,T_1-\epsilon)\times \p M_0),$$ such that the solution $u_{f_j}$ to \eqref{pf0} with boundary data $f_j$ satisfies: 
\begin{itemize}
\item[(i)]{$|u_{f_j}(p)-1|\leq \frac{C}{j}$ for all $j \in \N$, where $C$ is a constant that is independent of $j$.}
\item[(ii)]{$|\langle w,\nabla u_{f_j}(p)\rangle -{\rm i} \,j\,\langle w, \xi \rangle| \leq C$ for all $j \in \N$ and all $w \in T_pM$. Here ${\rm i}$ is the imaginary unit and the constant $C$ is independent of $j$.}
\end{itemize}
\end{lemma}
\begin{proof}[Proof of Theorem~\ref{thm_main}]
Let $p=(T_1,x_0)\in \mathbb D$, where $\mathbb D$ is as defined in \eqref{recovery_D} and let $U$ be a small neighborhood of $p$ in $\mathbb D$. Denote by $\pi:M\to \R$ the time projection given by $\pi(t,x)=t$. We fix $\epsilon>T_1-\inf_{q \in U} \pi(q)$ such that for each $v \in L^+U$, the null geodesic $\gamma_v$ is disjoint from $[T_1-\epsilon,T_1+\epsilon]\times \p M_0$. This can always be achieved for $U$ sufficiently small. 

Invoking Proposition~\ref{prop_final}, we observe that given each point $q\in U$ and each $f \in H^{\frac{n+1}{2}}_0((-T,T_1-\epsilon)\times \p M_0)$, there holds:
\begin{equation}\label{iden_U} u_f^{(1)}(q)= \omega(q)\,u_f^{(2)}(q), \quad \text{on $U$},\end{equation}
where $\omega(q)$ is independent of $f$.

We claim that $\omega \in C^{\infty}(U)$. To see this, we write $\widetilde{T_1}=T_1-\epsilon$ and let $q \in U$ be arbitrary. Let $\xi \in L_q^+M$ and note that $\gamma_{q,\xi}$ hits the timelike boundary $\Sigma$ at both of its end points. Noting that $\gamma_v$ is disjoint from $\Sigma$ for any $t \in [T_1-\epsilon,T_1+\epsilon]$ independent of $q \in U$, we can apply Lemma~\ref{asymp_value}. Thus there exists $f \in C^{\infty}_c((-T,\widetilde{T_1})\times \p M_0)$ such that $u_f^{(2)}\neq 0$ in a sufficiently small neighborhood $\widetilde{U}\subset U$ of the point $q$. Since $\omega= u_f^{(1)}(u_f^{(2)})^{-1}$ on $\widetilde{U}$, and since both expressions are smooth and the denominator is non-vanishing on $\widetilde{U}$, we conclude that the function $\omega$ is smooth near $q$. Since $q\in U$ was arbitrary, we conclude that $\omega \in C^{\infty}(U)$.  

Next, fixing $q\in U$, we let $f \in C^{\infty}_c((-T,\widetilde{T_1})\times \p M_0)$ be arbitrary. Noting that $u_f^{(1)}=\omega u_f^{(2)}$ on $U$, we write
\begin{equation}\label{final_eq}0=(\Box+V_1)u_h^{(1)}=  u_f^{(2)}\Box \omega -2 \langle \nabla \omega,\nabla u_f^{(2)}\rangle + (V_1-V_2)\,\omega\, u_f^{(2)}\quad \text{on $U$},\end{equation}
where we used the fact that $(\Box+V_j)u_f^{(j)}=0$, $j=1,2$. Let $\xi \in L^+_qM$ be an arbitrary future pointing null vector. Let $v=(q,\xi)$ and note again that the end points of the null geodesic $\gamma_v$ lie on $\Sigma$. Thus, in view of Lemma~\ref{asymp_value}, there exists a family of boundary data $\{f_j\}_{j\in \N}\subset C^{\infty}_c((-T,\widetilde{T_1})\times \p M_0)$ such that (i)--(ii) hold. Using these inequalities, substituting $f_j$ into \eqref{final_eq}, dividing by $j$ and taking the limit $j\to \infty$, we deduce that 
$$ \langle \nabla \omega(q), \xi\rangle =0.$$
Finally, since $\xi\in L_q^+M$ is an arbitrary future pointing null vector, we deduce that $\nabla \omega(q) =0$. As $q \in U$ is arbitrary, we conclude that $\omega$ is a non-zero constant function on $U$. That $\omega$ can not be identically zero follows from the fact that there is $f\in C^{\infty}_c((-T,\widetilde{T_1})\times \p M_0)$ such that $u_f^{(1)}(p)$ is non-zero, thanks to (i) in Lemma~\ref{asymp_value}. Thus, equation \eqref{final_eq} reduces to 
$$ (V_1-V_2) \,u_f^{(2)} = 0\quad \text{on $U$},$$
and therefore $V_1(p)=V_2(p)$. Here we are using (i) in Lemma~\ref{asymp_value} again to guarantee that there is $f\in C^{\infty}_c((-T,\widetilde{T_1})\times \p M_0)$ such that $u_f^{(2)}(p)$ is non-zero. Finally, since $p \in \mathbb D$ was arbitrary, we conclude that $V_1=V_2$ on the set $\mathbb D$.
\end{proof}

\medskip\paragraph{\bf Acknowledgements}
S.A  gratefully acknowledges support from NSERC grant 488916. A.F. was supported by EPSRC grant EP/P01593X/1. L.O.
was supported by EPSRC grants EP/R002207/1 and EP/P01593X/1.
\appendix

\section{}
\label{appendix A}
In this appendix, we prove that a smooth Lorentzian manifold $(\mathcal M,\textsl g)$ satisfying the assumption (H1) and with a timelike boundary, admits an isometric embedding 
$$\Phi:[-T,T]\times M_0\to \mathcal U\subset \mathcal M,$$
where $\mathcal U$ contains the set $\{q\in \mathcal M\,:\, \tau(q)\in [-T,T]\}$ and the metric $g=\Phi^*\textsl g$ is of the form \eqref{splitting}.

We start by considering $\tau:\mathcal M \to \R$, a proper, smooth, surjective function such that $\nabla \tau$ is timelike. Using \cite[Lemma 2.1]{HU} it follows that $\mathcal M$ is diffeomorphic to $\R \times M_0$ where $M_0=\{q \in\mathcal M\,:\,\tau(q)=0\}$. To conclude the proof we need to show the existence of a smooth function $\phi:\mathcal U\to \R$ such that $\nabla \phi$ is timelike and tangential to $\p \mathcal M$ for some $\mathcal U$ containing $\{q\in \mathcal M\,:\, \tau(q)\in [-T,T]\}$.

For $j=1,2,$ we define $\mathcal U_j= \{q\in \mathcal M\,:\,\tau(q) \in[-T-j,T+j]\}$ and $\Gamma_j=\mathcal{U}_j \cap \p \mathcal M$. We write $(z^0,z')$ for the boundary normal coordinates near $\Gamma_2$, so that 
$$ \Gamma_2= \{z^0=0\},\quad \text{and}\quad g(z^0,z')= (dz^0)^2+h(z^0,z'),$$
where $h(z^0,\cdot)$ is a family of smooth Lorentzian metrics on $\Gamma$ that smoothly depend on $z^0$. Then, define the function $\psi_0$ near $\Gamma_2$ in boundary normal coordinates by 
$$ \psi_0(z^0,z')= \tau(0,z'),$$  
and observe that $\nabla \psi_0$ is tangential and future pointing on $\Gamma_2$. Indeed, this suggests that the choice $\phi=\psi_0$ is suitable near $\Gamma_2$.

To complete our construction, we need to smoothly adjust $\psi_0$ as we move away from $\Gamma_2$ in such a way that $\nabla \phi$ remains timelike. There exists a nontangential future pointing timelike vector field $Z$ on $\Gamma_2$. For $\epsilon>0$ small, we define a map $F:[0,\epsilon]\times \Gamma_2 \to \mathcal M$ by 
$$ F(y^0,y')= \exp_{y'}(y^0 Z(y')).$$ 
It is straightforward to see that when $\epsilon>0$ is sufficiently small, $F$ is a diffeomorphism onto its image in $\mathcal M$. Moreover, since $Z$ is timelike and future pointing, it follows that $\nabla y^0$ is timelike and future pointing for small enough $\epsilon>0$.
(Here the gradient is computed with respect to $F^*\textsl g$.) Next, define a smooth function $\psi$ on $[0,\epsilon]\times \Gamma_2$ by
$$ \psi(y^0,y')= \chi(\epsilon^{-1}y^0)\left(\psi_0(y^0,y')-a\right) + (1- \chi(\epsilon^{-1}y^0))\tau(F(y^0,y')),$$
where $a$ is a positive number and $\chi:[0,\infty)\to \R$ is a smooth decreasing function satisfying $\chi(t)=1$ for $t< \frac{1}{4}$ and $\chi(t)=0$ for $t>\frac{1}{2}$.

We now define the desired smooth function $\phi$ on the set $\mathcal U=\mathcal U_1$ by
\begin{align*}
\phi &= \tau \quad\quad\quad\quad \text{on $\mathcal U_1 \setminus F([0,\epsilon]\times \Gamma_2)$},\\
\phi &= \psi\circ F^{-1} \quad \text{on $\mathcal U_1\cap F([0,\epsilon]\times \Gamma_2)$}.
\end{align*}
Since $\chi$ is decreasing and since $\nabla y^0$ and $\nabla \tau$ are timelike and future pointing, it follows that $\nabla \phi$ is timelike and future pointing on $\mathcal U_1$ for all large enough $a>0$. It is also straightforward to see that $\nabla\phi$ is tangential to the boundary. 
\section{}
\label{appendix B}
This appendix is concerned with proving Lemma~\ref{asymp_value}. We will in fact prove the lemma for quite general Lorentzian manifolds assuming only that $M$ is of the form \eqref{cylinder} and that the metric $g$ is of the form \eqref{splitting}. We begin by briefly recalling the classical Gaussian beam construction for the wave equation that was introduced in \cite{Ba,Ralston}.  The presentation here follows that of \cite{FO}. Here, by Gaussian beams  we refer to approximate solutions to the equation 
$$ (\Box+ V)u=0, \quad \text{on $M=[-T,T]\times M_0$},$$
that concentrate on an inextendible null geodesic $\gamma_v:[a,b]\to M$ where $v=(p,\xi) \in L^+M$ is such that the end point $\gamma_v(a)$ lies on the timelike boundary $\Sigma=(-T,T)\times \p M_0$. 

We begin by extending the manifold $M_0$ into a slightly larger manifold $\hat{M}_0$ and define $\hat{M}=(-T,T)\times \hat{M}_0$. We also extend the metric $g$ smoothly to $\hat{M}$ such that \eqref{splitting} holds over the extended manifold $\hat{M}$. Note that the null geodesic $\gamma_v$ also extends into $\hat{M}$ with end points on $(-T,T)\times \p \hat{M}_0$. We write $\gamma_v:[a',b']\to \hat{M}$ for this extended null geodesic. 

In order to recall the expression of Gaussian beams in local coordinates, we first briefly recall the well-known Fermi coordinates near a null geodesic. We refer the reader to \cite[Section 4.1, Lemma 1]{FO} for the proof.

\begin{lemma}[Fermi coordinates]
\label{fermi}
Let $\gamma_v: [a',b'] \to \hat{M}$ be a null geodesic on $\hat{M}$ parametrized as given by \eqref{affine} and whose end points lie on $(-T,T)\times \p \hat{M}_0$. Given each $\hat{a}\in (a',a)$ and $\hat{b}\in (b,b')$, there exists a coordinate neighborhood  $(U,\psi)$ of $\gamma([\hat{a},\hat{b}])$, with the coordinates denoted by $(y^0:=s,y^1,\ldots,y^n)=(s,y')$, such that:
\begin{itemize}
\item[(i)] {$\psi(U)=(\hat{a},\hat{b}) \times B(0,\delta')$ where $B(0,\delta')$ is the ball in $\mathbb{R}^{n}$ centered at the origin with a small radius $\delta' > 0$.}
\item[(ii)]{$\psi(\gamma(s))=(s,\underbrace{0,\ldots,0}_{n \hspace{1mm}\text{times}})$}. 
\end{itemize}
Moreover, the metric tensor $g$ written in this coordinate system satisfies  
    \begin{align}\label{g_on_gamma}
g|_\gamma = 2ds\otimes dy^1+ \sum_{\alpha=2}^n \,dy^\alpha\otimes \,dy^\alpha,
    \end{align}
and $\frac{\p}{\p y^i} g_{jk}|_\gamma = 0$ for $i,j,k=0,\ldots,n$. Here, $|_\gamma$ denotes the restriction on the curve $\gamma$.
\end{lemma}

In Fermi coordinate, Gaussian beams can be defined via the ansatz,
\begin{equation}\label{tau_pos} \mathcal U_\lambda(y) = e^{{\rm i}\lambda \phi(y)} A_\lambda(y) \quad \text{for} \quad \lambda>0\end{equation}
and
\begin{equation}\label{tau_neg} \mathcal U_\lambda(y) = e^{-{\rm i}\lambda\bar{\phi}(y)}\bar{A}_\lambda(y) \quad \text{for}\quad\lambda<0,\end{equation}
where the phase and amplitude functions $\phi$ and $A_\lambda$ are defined by
\begin{equation}\label{phase-amplitude}
\begin{aligned}
 \phi(s,y') = \sum_{j=0}^{N} \phi_j(s,y')& \quad \text{and} \quad A_\lambda(s,y')= \chi(\frac{|y'|}{\delta'}) \sum_{j=0}^{N} \lambda^{-j}a_{j}(s,y'),\\
&a_j(s,y')=\sum_{k=0}^{N} a_{j,k}(s,y').
\end{aligned}
\end{equation}
Here, given each $j,k=0,\ldots,N$, the term $\phi_j$ is a complex-valued homogeneous polynomial of degree $j$ in the variables $y^{1},\ldots, y^n$ and $a_{j,k}$ is a complex valued homogeneous polynomials of degree $k$ with respect to the variables $y^1,\ldots,y^n$. Finally, the function $\chi \in C^{\infty}_c(\R)$ is a non-negative function that satisfies $\chi(t)=1$ for $|t| \leq \frac{1}{4}$ and $\chi(t)=0$ for $|t|\geq \frac{1}{2}$.

The phase terms $\phi_j$ and the amplitudes $a_j$ with $j=0,1,2,\ldots,N$ are determined iteratively by solving ODEs along the null geodesic. This comes as a result of a WKB analysis for the conjugated wave operator in the semi-classical parameter $\lambda$, that enforces
\[
\begin{aligned}
&\frac{\p^{|\alpha|}}{\p y'^\alpha}\langle d\phi,d\phi\rangle=0 \quad \text{on $(\hat a,\hat b)\times\{y'=0\}$,}\\ 
&\frac{\p^{|\alpha|}}{\p y'^{\alpha}}\left( 2\langle d\phi,da_j \rangle- (\Box\phi)a_j + {\rm i} (\Box+V)a_{j-1}\right)=0\quad \text{on $(\hat a,\hat b)\times\{y'=0\}$},
\end{aligned}
\]
for all $j=0,1,\ldots,N$ and all multi-indices $\alpha=(\alpha_1,\ldots,\alpha_n) \in \{0,1,\ldots\}^n$ with $|\alpha|=\alpha_1+\ldots+\alpha_n\leq N$.

We do not proceed to solve these equations here as this can be found in all the works mentioned above, but instead summarize the main properties of Gaussian beams as follows:
\begin{enumerate}
\item[(1)] {$\phi(s,0)=0$. }
\item [(2)]{$\Im(\phi)(s,y') \geq C |y'|^2$ for all points $y \in (\hat a,\hat b)\times B(0,\delta')$.}
\item [(3)]{$\| (\Box+V) \mathcal U_\lambda \|_{H^k(M)} \lesssim |\lambda|^{-N'},$ where $N'=\frac{N+1}{2}+\frac{n}{4}-k-2$.}
\end{enumerate}
Here, $\Im$ stands for the imaginary part of a complex number. We will need explicit expressions for the phase terms $\phi_0$, $\phi_1$, $\phi_2$ and the principal amplitude $a_{0,0}$. The phase terms are given by the expressions
\begin{equation}
\label{principal_phase}
\begin{aligned}
\phi_0(s,y')=0,\quad \phi_1(s,y')= y^1,\quad \phi_2(s,y')=\sum_{j,k=1}^n H_{jk}(s)y^jy^k,
\end{aligned}
\end{equation}  
where the symmetric complex valued matrix $H$ solves the Riccati equation
\begin{equation}\label{riccati}
\frac{d}{ds} H + HCH + D=0, \quad \forall s \in (\hat a,\hat b), \quad H(0)=H_0,\quad \Im H_0>0.
\end{equation}
where $C$ and $D$ are the matrices defined through
\begin{equation}\label{Cmatrix}
\begin{cases}
C_{11}= 0&\\
C_{jj}=2& \quad j=2,\ldots,n, \\
 C_{jk}=0& \quad \text{otherwise,}
\end{cases}
\qquad \text{where $D_{jk}= \frac{1}{4} \frac{\p^2 g^{11}}{\p y^j \p y^k}$}.
\end{equation}
We recall the following result from \cite[Section 8]{KKL} regarding solvability of the Riccati equation.
\begin{lemma}
\label{ricA}
Let $H_0=Z_0Y_0^{-1}$ be a symmetric matrix with $\Im H_0 > 0$.
The Riccati equation (\ref{riccati}), together with the initial condition $H(0) = H_0$, admits a unique solution $H(s)$ for all $s \in [\hat a,\hat b]$. We have $\Im H(s)>0$ for all $s \in [\hat a,\hat b]$ and $H(s)=Z(s)Y^{-1}(s)$, where the matrix valued functions $Z(s)$, $Y(s)$ solve the first order linear system
$$ \frac{d}{ds} Y = CZ\quad \text{and}\quad  \frac{d}{ds} Z = -DY, \quad \text{subject to} \quad Y(0)=Y_0,\quad Z(0)=H_0.$$ 
Moreover, the matrix $Y(s)$ is non-degenerate on $[\hat a,\hat b]$, and there holds
$$
\det(\Im H(s)) \cdot |\det(Y(s))|^2=\det(\Im(H_0)).
$$
\end{lemma}
The principal part of the amplitude, that is the function $a_{0,0}$ is given by the expression:
\begin{equation}
\label{principal_amp}
a_{0,0}(s)=(\det Y(s))^{-\frac{1}{2}},
\end{equation}
where $Y(s)$ is as described above, by Lemma~\ref{ricA}. 

As for the remainder of the terms $\phi_j$ with $j\geq 3$ and the rest of the amplitude terms $a_{j,k}$ with $j,k$ not both zero, we recall from \cite{FO} that they solve first order ODEs along the null geodesic $\gamma$ and can be determined uniquely by fixing their initial values to be zero at the point $s=0$. This completes our review of Gaussian beams.

Next we assume as in the statement of Lemma~\ref{asymp_value} that $\epsilon>0$ is small enough so that $\gamma_v$ is disjoint from $[T_1-\epsilon,T_1+\epsilon]\times \p M_0$. We show that it is possible to choose $f \in C^{\infty}_c((-T,T_1-\epsilon)\times \p M_0)$, such that the solution $u$ to \eqref{pf0} with Dirichlet data $f$ is asymptotically close to the Gaussian beams $\mathcal U_\lambda$ on the subset $(-T,T_1+\frac{\epsilon}{2})\times M_0$. To this end, we consider for each $\lambda>0$, the Gaussian beam construction $\mathcal U_\lambda$ described above with 
\begin{equation}\label{N_order}
N=\lceil\frac{3n}{2} \rceil +10,
\end{equation}
where $\lceil x\rceil$ stands for the smallest integer that is greater than or equal to $x$ and we recall that the order $N$ is related to the Taylor series approximations of the phase and amplitude terms, see \eqref{phase-amplitude}. Let $\eta:\R\to \R$ be a smooth function such that $\eta(t) =1$ for $t<T_1+\frac{\epsilon}{2}$ and $\eta(t)=0$ for $t>T_1+\epsilon$. Next, let $u_\lambda$ be the solution to the wave equation
\begin{equation}\label{pf_gaussian}
\begin{aligned}
\begin{cases}
\Box u_\lambda+Vu_\lambda=0\,\quad &\text{on $M$},
\\
u_\lambda=f_\lambda=\eta\,\mathcal U_\lambda\,\quad &\text{on $\Sigma=(-T,T)\times \p M_0$,}\\
u_\lambda(-T,x)=\p_{t}u_\lambda(-T,x)=0 \,\quad &\text{on $M_0$,}
\end{cases}
    \end{aligned}
\end{equation}
Observe that given $\delta'$ sufficiently small, there holds 
$$\supp f_\lambda \subset (-T,T_1-\epsilon)\times \p M_0.$$ 
In view of property (3) in the construction of the Gaussian beam and the Sobolev embedding 
$$C^1((-T,T_1+\frac{\epsilon}{2})\times M_0)\subset H^{\frac{n+5}{2}}((-T,T_1+\frac{\epsilon}{2})\times M_0),$$
 we conclude that
\begin{equation}\label{error_bound}
\|u_\lambda - \mathcal U_\lambda\|_{C^1((-T,T_1+\frac{\epsilon}{2})\times M_0)}\leq \|(\Box+V)\mathcal U_\lambda\|_{H^{\frac{n+5}{2}}(M)}\leq \frac{C}{\lambda},
\end{equation} 
for some constant $C>0$ that is independent of $\lambda$. We are now ready to prove Lemma~\ref{asymp_value}.

\begin{proof}[Proof of Lemma~\ref{asymp_value}]

Note that by the hypothesis of the lemma, the end point $\gamma(a)$ of the inextendible null geodesic $\gamma_v:[a,b]\to M$ , lies on $\Sigma$. For $\lambda>0$, we consider the above Gaussian beams $\mathcal U_\lambda$ with $N$ satisfying \eqref{N_order} and $\delta'$ sufficiently small so that $\supp\, (\eta\,\mathcal U_\lambda)|_\Sigma \subset (-T,T_1-\epsilon)\times \p M_0$. Let 
$$f_\lambda=\eta\, \mathcal U_\lambda|_{\Sigma} \in C^{\infty}_c((-T,T_1-\epsilon)\times \p M_0),$$ 
and denote by $u_\lambda$, the unique solution to \eqref{pf0} subject to the Dirichlet data $f_\lambda$. In view of the explicit expressions \eqref{phase-amplitude}, \eqref{principal_phase} and \eqref{principal_amp}, there holds
$$\mathcal U_\lambda (p) = e^{{\rm i}\lambda \phi(0)}A_\lambda(0)=a_{0,0}(0)=(\det Y(0))^{-\frac{1}{2}},$$
and 
$$ \nabla \mathcal U_\lambda(p)={\rm i}\dot{\gamma}_v(0) A_\lambda(0)= {\rm i}\dot{\gamma}_v(0)  (\det Y(0))^{-\frac{1}{2}}.$$
The claim follows trivially from the latter two identities together with the error bound \eqref{error_bound}.
\end{proof}

\ifoptionfinal{}{
}

\begin{thebibliography}{99}

\bibitem{Alexakis}
S. Alexakis and A. Shao, Global uniqueness theorems for linear and nonlinear waves, J. Func. Anal. 269
(2015), no. 11, 3458–3499.

\bibitem{AB08} S. Alexander, R. Bishop, Lorentz and semi-Riemannian spaces with Alexandrov curvature
bounds, Comm. Anal. Geom. 16 (2008),251--282.

\bibitem{Al}
S. Alinhac, Non-unicit\'{e} du probl\'{e}me de Cauchy, Ann. of Math., 117 (2) (1983), 77--108.

\bibitem{AL04}
M. Anderson, A. Katsuda, Y. Kurylev, M. Lassas, M. Taylor, Boundary regularity for the Ricci equation, geometric convergence, and Gel’fand’s inverse boundary problem. Invent. math. 158, 261–321 (2004).

\bibitem{AH98}
L. Andersson, R. Howard, Comparison and rigidity theorems in semi-Riemannian geometry, Comm. Anal. Geom. 6 (1998), 819--877.

\bibitem{AK}
S. Alexander, W. Karr, spacetime convex functions and sectional curvature, Proceedings of the International Meeting on Lorentzian Geometry, M\'{a}laga, 2016

\bibitem{Ba}
V. Babich, V. Ulin, The complex spacetime ray method and quasi-photons, Zap. Nauch
Semin. LOMI 117 (1981), 5--12 (Russian)

\bibitem{BG}
N. Burq, P. G\'{e}rard, Condition n\'{e}cessaire et suffisante pour la controlabilit\'{e} exacte des ondes. (French)
[a necessary and sufficient condition for the exact controllability of the wave equation]. C. R. Acad.
Sci. Paris S\'{e}r. I Math. 325(7), 749–752 (1997)

\bibitem{BLR}
C. Bardos, G. Lebeau, and J. Rauch, Sharp sufficient conditions for the observation, control and stabilization of
waves from the boundary SIAM J. Contr. Opt. 30 1024--65 (1992).

\bibitem{BE1}
J. K. Beem, P. E. Ehrlich, Singularities, incompleteness and the
Lorentzian distance function, Math. Proc. Camb. Phil. Soc. 85, 161--178 (1979).

\bibitem{Beem}
J. K. Beem, P. E. Ehrlich, K. L. Easley, Global Lorentzian geometry (Second). Marcel Dekker, Inc., New York (1996).

\bibitem{Bel87}
M. Belishev, An approach to multidimensional inverse problems for the wave equation, Dokl. Akad. Nauk SSSR, 297 (1987), 524--527.

\bibitem{Bel}
M. Belishev, Recent progress in the boundary control method, Inverse Problems, 23 (2007), R1--R67.

\bibitem{Bel92} 
M. Belishev, Y. Kurylev, To the reconstruction of a Riemannian manifold via its spectral data (BC-method), Comm. Partial Differential Equations, 17 (1992), 767--804.

\bibitem{BS}
A.N. Bernal, M. S\'{a}nchez, Globally hyperbolic spacetimes can be defined as ‘causal’ instead of
‘strongly causal’, Class. Quant. Grav. 24, 745 (2007),

\bibitem{Blago}
A. S. Blagovestchenskii, “A one-dimensional inverse boundary value problem for
a second order hyperbolic equation”. Zap. Naučn. Sem. Leningrad. Otdel. Mat. Inst.
Steklov. (LOMI) 15, pp. 85–90 (1969).

\bibitem{DKSU}
D. Dos Santos Ferreira, C. E. Kenig, M. Salo, G. Uhlmann, Limiting Carleman weights and anisotropic inverse problems. Inventiones Mathematicae,178(1), 119--171, (2009). 

\bibitem{DKLS}
D. Dos Santos Ferreira, Y. Kurylev, M. Lassas, M. Salo, The Calder\'{o}n problem in transversally anisotropic geometries. J. Eur. Math. Soc. (JEMS), 18(11), 2579--2626, (2016). 


\bibitem{E1}
G. Eskin, Inverse hyperbolic problems with time-dependent coefficients, Commun. Partial Diff. Eqns., 32  (11) (2007), 1737--1758.

\bibitem{E2}
G. Eskin, Inverse problems for general second order hyperbolic equations with time-dependent coefficients,   Bull. Math. Sci., 	7 (2017), 247--307.

\bibitem{E3}
G. Eskin, Lectures on Linear Partial Differential Equations, Graduate Studies in Mathematics, vol. 123, AMS (2011)

\bibitem{Evans}
L. C. Evans, Partial Differential Equations, Volume 19 of Graduate studies in mathematics, American Mathematical Soc., 2010.

\bibitem{FIKO}
A. Feizmohammadi, J. Ilmavirta, Y. Kian, L. Oksanen, Recovery of time dependent coefficients from boundary data for hyperbolic equations, Journal of Spectral Theory, To appear (2020).

\bibitem{FIO}
A. Feizmohammadi, J. Ilmavirta, L. Oksanen, The light ray transform in stationary and static Lorentzian geometries,  J. Geom. Anal, To appear (2020).

\bibitem{FO}
A. Feizmohammadi, L. Oksanen, Recovery of zeroth order coefficients in non-linear wave equations, J. Inst. Math. Jussieu, To appear (2020).

\bibitem{GI}
G. Gibbons, A. Ishibashi, Convex functions and spacetime geometry, Classical Quantum Gravity 18 (2001), no. 21, 4607--4627.

\bibitem{HE}
S. W. Hawking, G. F. R. Ellis, “The Large Scale Structure Of Space-time,” Cambridge University Press (1973).

\bibitem{HU}
P. Hintz, G. Uhlmann, Reconstruction of Lorentzian Manifolds from Boundary Light
Observation Sets, International Mathematics Research Notices, Vol. 2019, No. 22, pp. 6949--
6987.

\bibitem{HUZ}
P. Hintz, G. Uhlmann, J. Zhai, An inverse boundary value problem for a semilinear wave equation on Lorentzian manifolds, arXiv preprint (2020).

\bibitem{Ho3}
L. H\"ormander, The Analysis of linear partial differential operators, Vol III, Springer-Verlag, Berlin, Heidelberg, 1983.

\bibitem{Ho4}
L. H\"ormander, The Analysis of linear partial differential operators, Vol IV, Springer-Verlag, Berlin, Heidelberg, 1983.

\bibitem{Ho}
L.  H\"ormander, A uniqueness theorem for second order hyperbolic differential equations, Comm.Partial Differential Equations, 16 (1991) 789--800.

\bibitem{I1}
V. Isakov. An inverse hyperbolic problem with many boundary measurements. Commun. Partial Differ.
Equ. 16, 1183--1195 (1991)

\bibitem{KKL} 
A. Katchalov, Y. Kurylev, M. Lassas, Inverse boundary spectral problems, Chapman \& Hall/CRC Monogr. Surv. Pure Appl. Math., 2001.


\bibitem{K}
R. Kulkarni, The values of sectional curvatures in indefinite metrics, Comment. Math. Helv. 54 (1979), 173--176.

\bibitem{KL2}
Y. Kurylev, M. Lassas. Hyperbolic inverse boundary value problems and time-continuation of the
non-stationary Dirichlet-to-Neumann map. Proc. R. Soc. Edinb. 132, 931--949 (2002)

\bibitem{KOP}
Y. Kurylev, L. Oksanen, G. P. Paternain, Inverse problems for the connection Laplacian, J. Differential Geom., 110 (2018), no. 3, 457--494. 

\bibitem{KLU}
Y. Kurylev, M. Lassas, G. Uhlmann, Inverse problems for Lorentzian manifolds and non-linear hyperbolic equations. Inventiones mathematicae, 212(3):781-857, 2018.

\bibitem{LLT}
I. Lasiecka, J-L. Lions, R. Triggiani, Non homogeneous boundary value problems for second order hyperbolic operators J. Math. Pures Appl., 65 (1986), 149--192.

\bibitem{Lassas}
M. Lassas, Inverse problems for linear and non-linear hyperbolic equations. Proc. Int. Cong. of Math. 2018, Rio de Janeiro, 3, 2018.


\bibitem{LOSU}
M. Lassas, L. Oksanen, P. Stefanov, and G. Uhlmann. The Light Ray Transform on Lorentzian Manifolds. Communications in Mathematical Physics (2020).

\bibitem{Oneill} 
O’Neill, B. (1983). Semi-Riemannian geometry. New York: Academic Press Inc.

\bibitem{Ralston}
J. Ralston, Gaussian beams and the propagation of singularities. Studies in Partial Differential Equations, MAA Studies in Mathematics 23, 206--248 (1983)

\bibitem{RS}
A. G. Ramm and J. Sj\"ostrand, An inverse problem of the wave equation, Math. Z., 206 (1991), 119--130.

\bibitem{Ring}
H. Ringstr\"om, The Cauchy Problem in General Relativity, ESI Lectures in Mathematics
and Physics. Zurich: European Mathematical Society Publishing House (2009).

\bibitem{Rob}
L. Robbiano,  Th\'{e}or\'{e}me  d'unicit\'{e}  adapte  au  controle  des  solutions  des  probl\'{e}mes  hyperboliques, Comm. Partial Differential Equations, 17 (1992) 699--714.

\bibitem{RZ}
L. Robbiano, C. Zuily. Uniqueness in the Cauchy problem for operators with partially holomorphic
coefficients. Invent. Math. 131(3), 493--539 (1998)

\bibitem{Sal}
R. Salazar. Determination of time-dependent coefficients for a hyperbolic inverse problem. Inverse
Probl. 29(9), 095015 (2013)

\bibitem{Shao}
A. Shao, On Carleman and observability estimates for wave equations on time-dependent domains,
Proc. Lond. Math. Soc. 119 (2019), no. 4, 998--1064.

\bibitem{St}
P. Stefanov, Uniqueness of the multi-dimensional inverse scattering problem for time-dependent potentials, Math. Z., 201 (4) (1989), 541--559.

\bibitem{Stefanov1}
P. Stefanov, Support theorems for the light ray transform on analytic Lorentzian manifolds, Proc. Amer. Math. Soc., 145, pp. 1259--1274, 2017.

\bibitem{SY}
P. Stefanov, Y. Yang, The inverse problem for the Dirichlet-to-Neumann map on Lorentzian manifolds. Anal. PDE. 2018;11(6):1381--1414.

\bibitem{Tataru}
D. Tataru, Unique continuation for solutions to PDE; between H\"{o}rmander’s theorem and Holmgren’s theorem,
Commun. Partial Diff. Eqns., 20 (1995), 855--884.

\end{thebibliography}
\end{document}